\documentclass[10pt]{article}
\usepackage{latexsym}
\usepackage[english]{babel}
\usepackage{graphicx,verbatim}
\usepackage{xcolor}
\usepackage{epstopdf}
\usepackage{amssymb}
\usepackage{algorithm}
\usepackage{algpseudocode}
\usepackage{amsmath,amsthm,enumerate,rotate}
\usepackage{hyperref}
\usepackage{tabularx,multirow,booktabs}

\def\be#1\ee{\begin{equation}#1\end{equation}}

\newtheorem{remark}{\bf Remark}[section]
\newcommand{\bq}{\begin{equation}}
\newcommand{\eq}{\end{equation}}

\def\bqa{\begin{eqnarray}}
\def\eqa{\end{eqnarray}}

\newcommand{\bd}{\begin{displaymath}}
\newcommand{\ed}{\end{displaymath}}
\newcommand{\ba}{\begin{eqnarray}}
\newcommand{\ea}{\end{eqnarray}}

\newcommand{\first}[1]{\textcolor{black}{#1}}
\newcommand{\second}[1]{\textcolor{black}{#1}}
\begin{document}

%\title{Kinetic models for opinion-driven epidemic dynamics on static networks}

%\title{Kinetic models for opinion-driven epidemic dynamics on networks}

\title{A second-order numerical scheme for optimal control of non-linear Fokker--Planck equations and applications in social dynamics}

\author{G. Albi\thanks{Department of Computer Science,
		University of Verona, Strada le Grazie 15, 37134 Verona, Italy. (giacomo.albi@univr.it)}\and E. Calzola \thanks{Department of Computer Science,
		University of Verona, Strada le Grazie 15, 37134 Verona, Italy. (elisa.calzola@univr.it)}}

%%%% Subject entries to be placed here %%%%
%\subject{applied mathematics, social sciences, partial differential equations, numerical
%simulations}

%%%% Keyword entries to be placed here %%%%
%\keywords{multi-agent systems, Boltzmann equation, }

\maketitle

\abstract{%
	%numerical methods for optiml control nonlocal fokker planck and applications to socio-economics dynamics
	In this work, we present a second-order numerical scheme to address the solution of optimal control problems constrained by the evolution of nonlinear Fokker--Planck equations arising from socio--economic dynamics. In order to design an appropriate numerical scheme for control realization, a coupled forward-backward system is derived based on the associated optimality conditions. The forward equation, corresponding to the Fokker--Planck dynamics, is discretized using a structure-preserving scheme able to capture steady states. On the other hand, the backward equation, modeled as a Hamilton--Jacobi--Bellman problem, is solved via a semi-Lagrangian scheme that supports large time steps while preserving stability. Coupling between the forward and backward problems is achieved through a gradient descent optimization strategy, ensuring convergence to the optimal control. Numerical experiments demonstrate second-order accuracy, computational efficiency, and effectiveness in controlling different examples across various scenarios in social dynamics. This approach provides a reliable computational tool for the study of opinion manipulation and consensus formation in socially structured systems.
}

\section{Introduction}
\second{In this manuscript we focus on the development of numerical methods for optimal control problems constrained by the evolution of nonlinear Fokker--Planck equations 
	% of the form
	% \begin{equation}\label{eq:protoFP}
		% \partial_t f =- \nabla \cdot\left(\mathcal{G}[f;\mathbf{u}](\mathbf{v},t) f - \nabla\cdot\left(\mathcal{D}(\mathbf{v})f\right)\right), 
		% \end{equation}
	% where $\mathbf{v}\in\Omega\subseteq \mathbb{R}^d$, $f(\mathbf{v},t)\geq 0$ and $t\geq 0$, and 
	% where the control $\mathbf{u}(\mathbf{v},t)\in \mathbb{R}^d$ is solution  to the following objective functional
	% \begin{equation}\label{eq:protocost}
		% \mathbf{u}=\arg\min_{\mathbf{u} \in U} \mathcal{J}(\mathbf{u};f_0),
		% \end{equation}
	% \second{with $\mathcal{J}$ as defined in \eqref{eq:cost_j}.} 
	% In particular, we aim at describing dynamics where $f(\mathbf{v},t)$ represents a continuous population, the controlled \second{nonlinear} drift $\mathcal G [\cdot]$ encodes non-local interactions among the agents' population, \first{with the notation $[\cdot]$ meaning functional dependence,} and the \second{nonlinear} function $\mathcal{D}(\cdot)\geq 0$, \first{$\mathcal{D}: \mathbb{R}^d \to \mathbb{R}_+^{d\times d}$ that} weights the influence of the diffusion.  
	% Fokker--Planck models of type \eqref{eq:protoFP} 
	arising in several fields of applications such as in plasma physics \cite{bobylev2000theory}, in biology \cite{carrillo2015non,degond2020fokker}, in socio-economic dynamics \cite{pinasco2017modeling,pareschi2014wealth,during2009boltzmann}, as well as optimization and learning methods \cite{carrillo2018analytical,fornasier20}.} The study of the nonlinear optimization problem constrained by PDEs %\eqref{eq:protoFP}-\eqref{eq:protocost} 
is motivated by the need to design efficient control strategies in these different areas of application domains. In particular in the last decade large interest has been shown in the context of control of multi-agent systems at various scales, see for example \cite{albi2025instantaneous,burger21,franceschi2024optimal,gugat2024turnpike,albi17,fornasier2014mean,during2009boltzmann}. In particular, the optimal control of Fokker--Planck models has been studied in several directions, for example in \cite{aronna2021,fleig2017optimal} from a theoretical perspective, or in a stochastic setting in \cite{breitenbach2020pontryagin,annunziato13}, as well as stabilization problems \cite{breiten2018control,bicego2024computation,alla2023state}. In this framework, the numerical realization of the control requires careful treatment because of nonlinearities involved that must be treated efficiently to achieve solutions with proper accuracy and within a reasonable computational time, see for example \cite{albi2024control,albi17,albi24_expint,aduamoah2022}. Here we address the control synthesis by means of reduced gradient method, solving the associated optimality conditions. This approach results in solving numerically the coupled system of the Fokker--Planck model and the associated backward equation, see for instance \cite{albi17}.

Specifically, the forward problem, corresponding to the evolution of the Fokker--Planck equation, is approximated with a Chang--Cooper scheme \cite{chang1970practical,mohammadi2015analysis}, where in this case for the nonlinear drift we follow the approach proposed in \cite{pareschi18}. 
The backward equation is modeled as a Hamilton--Jacobi--Bellman problem, for which we use a semi-Lagrangian scheme that supports large time steps while preserving stability \cite{falcone02,calzola23,alla2015efficient}.
Unlike standard Eulerian approaches, semi-Lagrangian schemes naturally handle the presence of advection and diffusion terms. They are explicit and stable without requiring any parabolic CFL condition, thus proving to be flexible tools for complex problems. These methods can achieve high-order accuracy by choosing appropriate interpolation techniques and, in the viscous case, by implementing different treatments of the diffusion terms \cite{ferretti10}. Moreover, their implementation adapts well also to various domain geometries and can be made efficient on unstructured meshes \cite{cacaceferretti}.
In order to optimize the associated cost, the forward and backward equations are coupled through an iterative reduced-gradient descent procedure, to provide the optimal control solution. Similar approaches for mean-field optimal control problem 
%such as \eqref{eq:protoFP}--\eqref{eq:protocost} 
have been studied in \cite{albi17,aduamoah2022,breiten2025optimal}, \second{here the main novelties of this manuscript consist in the combination of two second order schemes for the coupled forward and backward system in the multi-dimensional case, allowing increase of accuracy and efficiency, as well providing different numerical experiments in the context of collective dynamics in particular for opinions dynamics \cite{albi2014boltzmann,bicego2024computation,during2009boltzmann,toscani06,albi24,pinasco2017modeling}}.

%\newpage
%\cite{albi24} here we have the opinion dynamics model we want to control.
%\cite{albi16} here the aim is to drive opinions towards a given opinion target, but the process in of an evolving network. Particle level optimization discretizing in time without using lagrange multipliers.
%\cite{albi24_expint} here we have the numerical approximations of similar problems but with constant diffusion and 1-d, no contacts.
%\cite{annunziato13} here there is the reformulation from stochastic optimal control using SDEs to deterministic optimal control using the FP equation. Notice that there are some theoretical results with space-time dependent drift and diffusion functions
%\cite{butt24} here they solve the problem with constant diffusion and drift which is not nonlocal, using the Lagrange multipliers approach. CC is used for forward and backward
%\cite{falcone02} here there is the interpretation of sl schemes as Hopf-Lax-Oleĭnik representation formula for the exact viscosity solution of first order evolutive Hamilton--Jacobi equations. We instead have second-order formulas
%\cite{fleming06} classical book with theory of optimal control, especially stochastic optimal control
%\cite{lions71} a possible classical reference for the Lagrange multipliers approach in optimal control problems governed by PDEs
%\cite{roy18} 

The manuscript is structured as follows. Section \ref{sec:ocopi} presents the  PDE constrained optimization problem, and shows how to derive a forward-backward system of coupled equations and the optimality conditions for the control function. \second{Section \ref{sec:discretization} presents a gradient-descent approach to approximate the optimality conditions and introduces the combination of numerical schemes that implement it}. In particular, Section \ref{sec:secondordersl} introduces the \second{second-order} semi-Lagrangian discretization for the Hamilton--Jacobi equation, while Section \ref{sec:spscheme} presents the \second{second-order} \second{structure-preserving} scheme for the Fokker--Planck equation. \second{Neither scheme requires} a parabolic CFL condition for stability, allowing larger time steps, which is crucial when an unknown number of iterations of the gradient descent approach \second{is} necessary to reach the optimality condition. Section \ref{sec:numerics} presents some numerical simulations, in particular Section \ref{sec:num1} shows the order two of convergence of the \second{structure-preserving} scheme to the stationary solution of the Fokker--Planck equation, while Section \ref{sec:num2} is dedicated to showing the order two of convergence of the \second{two coupled} methods for the forward-backward system, given a fixed suboptimal control. Finally, Section \ref{sec:qualit} presents applications in social dynamics in two spatial dimensions: Section \ref{sec:opic} presents a model in which the first variable is the opinion on a given topic and the second one is the number of social media contacts, modeling the controlled process of opinion evolution on an evolving social network, while Section \ref{sec:dueopi} is dedicated to a \second{two-opinion} model, where the controlled evolution of two opinions (that influence each other) is studied. Section \ref{sec:conclusions} presents some concluding remarks and possible future research directions.

\section{Optimal control of non-local Fokker--Planck equations}\label{sec:ocopi}
Let $\mathbf{v}=(v_1,\dots,v_d)\in\Omega\subset\mathbb{R}^d$.
Given $T>0$, we will use the notation $\Omega_T = \Omega\times(0,T]$ for the spatiotemporal domain and $\vec{n}(\mathbf{v})$ to indicate the outward pointing normal vector in any point $\mathbf{v}\in\partial\Omega$. We want to numerically solve the mean-field optimal control problem 
\begin{equation}\label{eq:minimization}
	\min _{\mathbf{u} \in U} \mathcal{J}(\mathbf{u};f_0)
\end{equation}
\second{constrained to the evolution of the Fokker-Planck equation} 
\begin{equation}\label{eq:opifp}
	\begin{cases}
		\partial_t f =- \nabla \cdot\left(\mathcal{G}[f;\mathbf{u}](\mathbf{v},t) f - \nabla\cdot\left(\mathcal{D}(\mathbf{v})f\right)\right), \quad &(\mathbf{v},t) \in \Omega_T\\
		f(\mathbf{v},0) = f_0(\mathbf{v}), &\mathbf{v} \in \Omega\\
		\left(-\mathcal{G}[f,\mathbf{u}](\mathbf{v},t)f + \nabla\left(\mathcal{D}(\mathbf{v})f\right)\right) \cdot \vec{n}(\mathbf{v}) = 0, &\mathbf{v} \in \partial\Omega_,
	\end{cases}
\end{equation}
where \second{$f_0(\mathbf{v})\geq 0$ for all $\mathbf{v} \in \Omega$ such that $f_0\in L^2(\Omega)$ is a probability density function.} 
\second{The drift $\mathcal{G}$ is assumed to have the following structure}
\begin{equation}\label{eq:drift}
	\mathcal{G}[f;\mathbf{u}](\mathbf{v},t) = \boldsymbol{\mathcal{P}}[f](\mathbf{v},t) + \mathbf{u}(\mathbf{v},t),
\end{equation}
where $\boldsymbol{\mathcal{P}}[f](\mathbf{v},t)$ is a vector whose components are the operators $\mathcal{P}_i[f]$, that can be expressed as the sum of a local and a nonlocal term as
\begin{equation}\label{eq:drift_general}
	\mathcal{P}_i[f](\mathbf{v},t) = h_i(\mathbf{v})
	+\int_{\Omega} P_i(\mathbf{v},\mathbf{v}^*)(v_i^* - v_i)f(\mathbf{v}^*,t)\mathrm{d}\mathbf{v}^* , 
\end{equation}
with $h_i\in L^\infty(\Omega)$ and $P_i\in L^\infty(\Omega\times\Omega)$ for all $i=1,\dots,d$. The drift include an additive control term $\mathbf{u}=(u_1,\ldots,u_d)\in U$, where $U$ denotes the set of admissible controls defined as follows
\[
U := \left\{\mathbf{u} \in L^2(0,T;L^\infty(\Omega)) : \|\mathbf{u}\|_{L^2(0,T;L^\infty(\Omega))}\leq M\right\}
\]
for a fixed $M>0$.

The diffusion term in \eqref{eq:opifp} is characterized by the matrix
\begin{equation}\label{eq:diffusion}
	\mathcal{D}(\mathbf{v}) = \text{diag}\left(\mathcal{D}_1(\mathbf{v}), \dots, \mathcal{D}_d(\mathbf{v})\right),
\end{equation}
\first{where the functions $\mathcal{D}_i$ are such that for all $i=1, \dots, d$ $\mathcal{D}_i : \Omega \to \mathbb{R}_+$, $\mathcal{D}_i(\cdot)\in W^{1,\infty}(\Omega)$ and $\theta\leq \mathcal{D}_i(x)\leq \Theta$ for all $x\in\Omega$ and $0<\theta\leq \Theta<\infty$.}

The cost functional $\mathcal{J}$ is defined as
\begin{equation}\label{eq:cost_j}
	\mathcal{J}(\mathbf{u};f_0) = \frac{1}{2} \int_0^T\int_{\Omega} \left(\sum_{i=1}^d\lvert v_i-\overline v_i\rvert^2s_i(\mathbf{v}) + \gamma \lVert \mathbf{u}\rVert^2\right)f(\mathbf{v},t) \mathrm{d}\mathbf{v}\,\mathrm{d}t,
\end{equation}
where $\lVert \cdot \rVert$ indicates the Euclidean norm \second{in $\mathbb{R}^d$, and $s_i\in L^\infty(\Omega)$ a penalization function for $i=1,\dots,d$}.  
\begin{remark}
	\second{Within this framework of assumptions, the existence and the uniqueness of weak solution are guaranteed for the Fokker-Planck equation \eqref{eq:opifp} with $f\in L^2(0,T;H^1(\Omega))$}. \second{ Furthermore, it can be shown that starting from an initial datum $f_0$, that is a probability density, the zero-flux boundary conditions preserve the mass and the positivity of the solution, meaning that} $f(\mathbf{v},t)$ must be understood as a probability density, which means that
	\[
	\int_{\Omega}f(\mathbf{v},t)\mathrm{d}\mathbf{v} = 1 \quad \mbox{for all } t\in [0,T] \second{\quad \mbox{and}\quad f(\mathbf{v},t)\geq 0 \quad \mbox{for all} \quad(\mathbf{v},t)\in \Omega_T}.
	\]
	\second{Finally, existence of controls for \eqref{eq:minimization}--\eqref{eq:opifp} can be proven thanks to the structure of the objective functional \eqref{eq:minimization}, see for example \cite{albi17}.} 
	
	\second{We refer to \cite{korner2022second,aronna2021} for further details on the regularity and stability of such optimal control problem, as well as to \cite{frank2005nonlinear} for a general reference on the properties of non-linear Fokker-Planck equations.} 
\end{remark}

\subsection{Optimality conditions}
In order to solve the optimization problem \eqref{eq:minimization}-\eqref{eq:opifp} we resort to the \second{Lagrangian formalism (\cite{albi17,burger21,lions71}), to this end we introduce first the adjoint function $\psi(\cdot)$, such that $\psi\in L^2(0,T;H^1(\Omega))\cap C([0,T];L^2(\Omega))$ and $\partial_t\psi\in L^2(0,T;H^{-1}(\Omega))$}, and we introduce the Lagrangian as follows
\begin{equation}\label{eq:lagrangian}
	\begin{aligned}
		&\mathcal{L}(f,\mathbf{u},\psi) = \mathcal{J}(\mathbf{u};f_0)  \\
		&\quad-\int_0^T\int_{\Omega} \psi\left( \partial_t f + \nabla \cdot\left(\mathcal{G}[f;\mathbf{u}](\mathbf{v},t) f - \nabla\cdot\left(\mathcal{D}(\mathbf{v})f\right)\right)\right) \mathrm{d}\mathbf{v}\mathrm{d}t.
	\end{aligned}
\end{equation}

\second{Under the current assumptions, and the no-flux hypothesis, the sate-to-control map $u\to f(u)$ is continuously Fr\'echet differentiable on $u\in U$, with $f(u)$ the weak solution of \eqref{eq:opifp}, see \cite{aronna2021,korner2022second}
	for further details on this aspect}. Hence, we derive the optimality conditions by computing the Fr\'echet derivatives of the Lagrangian, meaning that we want to numerically solve
\begin{equation}\label{eq:frechetder}
	\begin{cases}
		D_f\mathcal{L}(f,\mathbf{u},\psi) = 0, \\
		D_{\mathbf{u}}\mathcal{L}(f,\mathbf{u},\psi) = 0, \\
		%D_{u_2}\mathcal{L}(f,u,\psi) = 0, \\
		D_\psi\mathcal{L}(f,\mathbf{u},\psi) = 0.
	\end{cases}
\end{equation}
The first equation in \eqref{eq:frechetder} results in the adjoint backward equation
\begin{equation}\label{eq:hj}
	\begin{split}
		-\partial_t \psi =\, & \mathcal{G}[f;\mathbf{u}](\mathbf{v},t)\cdot\nabla\psi + \left(\mathcal{D}(\mathbf{v})\nabla^\top\right)\cdot\nabla\psi \\&+ \sum_{i=1}^d\mathcal{Q}_i[f,\psi](\mathbf{v},t) 
		+ \frac{1}{2}\left(\sum_{i=1}^d\lvert v_i-\overline v_i\rvert^2s_i(\mathbf{v}) + \gamma \lVert \mathbf{u}\rVert^2 \right),
	\end{split}
\end{equation}
\first{with boundary condition 
	\[
	\mathcal{D}(\mathbf{v})\nabla\psi\cdot\vec{n}(\mathbf{v}) =0,
	\]
}
\first{coupled with the terminal condition $\psi (\mathbf{v},T) = 0$}, and where the nonlocal operators $\mathcal{Q}_i[f,\psi]$, $i=1,\dots,n$, are defined as
\[
\mathcal{Q}_i[f,\psi](\mathbf{v},t) = \int_{\Omega} P_i(\mathbf{v}^*,\mathbf{v})(v_i %\\
- v_i^*)f(\mathbf{v}^*,t)\partial_{v_i} \psi(\mathbf{v}^*,t)\mathrm{d}\mathbf{v}^*,
\]
\first{where $\partial_{v_i}$ indicates the $i$-th component of the nabla operator.}
The second equation gives the optimality condition
\begin{equation}\label{eq:optimcond}
	\gamma \mathbf{u}(\mathbf{v},t) + \nabla \psi(\mathbf{v},t) = 0
\end{equation}
and the third one corresponds to the Fokker--Planck equation \eqref{eq:opifp}.

In the following section we will detail the derivation of the numerical schemes involved, as well as their properties for the current optimal control problem.

\section{Discretization of the optimality system}\label{sec:discretization}
In order to find a candidate solution for the optimality condition system \eqref{eq:frechetder} we will proceed numerically using a reduced-gradient method. Hence, starting from $\mathbf{u}^{(0)}=\mathbf{0}$ and an initial guess for the density $f=f^{(0)}$, we aim to approach iteratively the solution of \eqref{eq:minimization}-\eqref{eq:opifp}.
At each iteration $k$\second{,} we approximate the solution to the forward problem \eqref{eq:opifp} with a semi-implicit \second{structure-preserving} scheme, in order to obtain $f^{(k)}$ over the whole time frame $[0,T]$. Subsequently, $f^{(k)}$ is passed \second{to} the backward equation \eqref{eq:hj}, where its solution $\psi^{(k)}$ is discretized with a semi-Lagrangian scheme. Then we perform a gradient descent step to update the control variable
\begin{equation}\label{eq:gradient}
	\mathbf{u}^{(k)} = \mathbf{u}^{(k-1)} -\second{\lambda^{(k)}} D_\mathbf{u}\first{\mathcal{L}}(f^{(k-1)},\mathbf{u}^{(k-1)},\psi^{(k)}),
\end{equation}
\second{where $\lambda^{(k)}$ is computed at every step using the Barzilai-Borwein method \cite{barzilai1988two}} and we iterate the process until \second{the difference $\lvert \mathcal{J}(\mathbf{u}^{(k)};f_0) - \mathcal{J}(\mathbf{u}^{(k-1)};f_0)\rvert<10^{-5}$}.

In the sequel, we will focus on the two-dimensional case in order to make the notation easier to read, but everything can be extended to a general $d$-dimensional framework. Let $\Omega = \Omega_1\times\Omega_2$, with $\Omega_1$ and $\Omega_2$ compact connected subsets of $\mathbb{R}$. We choose two discretization steps, $\Delta v_1,\Delta v_2 >0$, and since we can rewrite $\Omega_1 = [a_1,b_1]$ and $\Omega_2 = [a_2, b_2]$ we define the grid 
\begin{equation}\label{eq:grid}
	\begin{aligned}
		G_{\Delta v_1,\Delta v_2} = \{&\mathbf{v}_{i,j} = (v_{1,i},v_{2,j}) \mbox{ s.t.  for  } k=1,2\cr
		&\quad v_{k,i} = a_1+\Delta v_k(i+1/2),\,i=0,\dots,N_{v_k}-1\},
	\end{aligned}
\end{equation}
where $N_{v_1} = \second{\lceil}\lvert\Omega_1\rvert/\Delta v_1\second{\rceil\in\mathbb{N}_0}$ and $N_{v_2} = \second{\lceil}\lvert \Omega_2\rvert/\Delta v_2\second{\rceil \in \mathbb{N}_0}$. We also fix a time discretization step $\Delta t>0$, a number of time steps $N_T = \lceil T/\Delta t\rceil$, such that $t^n = n\Delta t$, $n=0,\dots,N_T$. Concerning the choice of $\Delta t$, both the semi-Lagrangian scheme and the semi-implicit \second{structure-preserving} one do not require a parabolic CFL condition in order to ensure stability, meaning that we can set the same time step for both methods equal to $\Delta t = \Delta v_1 \Delta v_2/(\second{\kappa}(\Delta v_1 + \Delta v_2))$, \first{with $\second{\kappa} = 2\max_{i,j}\lvert-\mathcal{G}[f;\mathbf{u}](\mathbf{v}_{i,j}) + \nabla \mathcal{D}(\mathbf{v}_{i,j})\rvert$, so} that positivity of the solution of the scheme for \eqref{eq:opifp} is preserved \first{(extension of \cite[Proposition 3]{pareschi18} in two dimensions)}. We will use the notation $\psi^n_{i,j}$ and $f^n_{ij}$ to indicate the numerical approximation of $\psi$ and $f$ in $\mathbf{v}_{i,j}\in G_{\Delta v_1, \Delta v_2}$ at time $t^n$, and we will indicate with $\pmb{\psi}^n=(\psi^n_{i,j})_{i,j}$ and $\mathbf{f}^n = (f^n_{i,j})_{i,j}$ the set of all their values on $G_{\Delta v_1, \Delta v_2}$.

\begin{remark}
	\first{The objective function is generally non-convex, due to the non-locality and the density-dependent drift present in the Fokker-Planck \eqref{eq:opifp}. Consequently, standard reduced gradient schemes, such as \eqref{eq:gradient}, converge only at first-order stationary points from arbitrary initializations. Hence, they do not guarantee convergence to a global minimizer, and different initial assumptions can lead to different local optima. Nevertheless, such adjoint-based scheme remains highly effective in practice, since a single gradient evaluation requires only two PDE solvers. Moreover, the use of Barzilai-Borwein step sizes, optionally combined with simple non-monotonic backtracking, produces inexpensive iterations that typically allow for rapid decay, \cite{barzilai1988two}.}
\end{remark}

\subsection{The second-order semi-Lagrangian scheme for the backward equation}\label{sec:secondordersl}
Considering $f, \mathbf{u}$ two given functions, we can write the second-order in time and space semi-Lagrangian approximation of \eqref{eq:hj} as in, for example, \cite{bonaventura21}. 
On the uniform grid $G_{\Delta v_1, \Delta v_2}$ it is possible to define $\mathbb{I}^r$, a piecewise polynomial interpolation operator of order $r$  constructed on the values of a grid function $\mathbf{g}$ such that  $\mathbb{I}^r[\mathbf{g}](\mathbf{v}_{i,j}) = g_{i,j}$. Let $K = [v_{1,i},v_{1,i+1}] \times [v_{2,j},v_{2,j+1}]$ be an element of the discretization defined in \eqref{eq:grid}, then for all $g\in C^{r+1}(K)$, $g(\mathbf{v}_{i,j}) = g_{i,j}$, the following local error estimate holds 
\[
\lVert g - \mathbb{I}^r[\textbf{g}]\rVert_{L^2(K)} \leq C\left(\max\{\Delta v_1, \Delta v_2\}\right)^{r+1},
\]
for a suitable constant $C>0$, leading to the global estimate
\[
\lVert g - \mathbb{I}^r[\textbf{g}]\rVert_{L^2(\Omega_1\times\Omega_2)} \leq C\left(\max\{\Delta v_1, \Delta v_2\}\right)^{r+1},
\]
if we also have $g\in H^{r+1}(\Omega)$
(we refer to \cite[Chapter 4.5]{quarteroni14} for a proof). In our simulations, we set $r=3$, in order to keep the interpolation error lower than the error due to the approximation of the stochastic characteristics. \second{In fact, since we can take $\Delta t = O(\Delta v_i)$, choosing $r=3$ the global interpolation error is of order $O(\Delta x^3) = O(\Delta t^3)$. In what follows, we use a method which is globally of order $\Delta t^2$ to approximate the characteristic curves, and since $\Delta t <1$ such error is the dominating one, leading to an overall order $2$ in time.}

Defining $\mathbf{V}=(V_1,V_2)$, for $s\in[t,t+\Delta t]$, the characteristic curves for equation \eqref{eq:hj} solve the following SDE
\begin{equation}\label{eq:charsde}
	\begin{cases}
		\mathrm{d}\mathbf{V} = \mathcal{G}[f;\mathbf{u}](\mathbf{V}(s),s)\mathrm{d}s + \sqrt{2\mathcal{D}(\mathbf{V}(s))}\mathrm{d}\mathbf{W}(s), \\
		\mathbf{V}(t+\Delta t) = \mathbf{v},
	\end{cases}
\end{equation}
where $\mathbf{W}$ is a standard two-dimensional Brownian motion. At each time step $t^n$ the second-order numerical approximation of the characteristic curves \eqref{eq:charsde} is given by $(y_{i,j}^{n,\ell})$, $\ell = 1,\dots,9$, where
\begin{equation}\label{eq:numchar}
	\begin{aligned}
		y_{i,j}^{n,\ell} = \mathbf{v}_{i,j} &+ \frac{\Delta t}{2}\left(\mathcal{G}[\mathbf{f}^n;\mathbf{u}^n](\mathbf{v}_{i,j},t^n) + \mathbb{I}^r[\first{\mathcal{G}}[\mathbf{f}^{n+1};\mathbf{u}^{n+1}]](y_{i,j}^{n,\ell})\right) \\
		&+ \sqrt{2\Delta t \mathcal{D}(\mathbf{v}_{i,j})}\mathbf{e}^\ell,
	\end{aligned}
\end{equation}
and the vectors $\mathbf{e}^\ell$, $\ell = 1,\dots,9$ are defined as
\[
\begin{aligned}
	&\mathbf{e}^1 = \begin{bmatrix}
		0 \\ 0
	\end{bmatrix}, \quad & \mathbf{e}^2 = \begin{bmatrix}
		+\sqrt{3} \\ 0
	\end{bmatrix}, \quad & \mathbf{e}^3 = \begin{bmatrix}
		-\sqrt{3} \\ 0
	\end{bmatrix}, \\
	&\mathbf{e}^4 = \begin{bmatrix}
		0 \\ +\sqrt{3}
	\end{bmatrix}, \quad & \mathbf{e}^5 = \begin{bmatrix}
		0 \\ -\sqrt{3} 
	\end{bmatrix}, \quad & \mathbf{e}^6 = \begin{bmatrix}
		+\sqrt{3} \\ +\sqrt{3}
	\end{bmatrix}, \\
	&\mathbf{e}^7 = \begin{bmatrix}
		-\sqrt{3} \\ +\sqrt{3}
	\end{bmatrix}, \quad & \mathbf{e}^8 = \begin{bmatrix}
		+\sqrt{3} \\ -\sqrt{3}
	\end{bmatrix}, \quad & \mathbf{e}^9 = \begin{bmatrix}
		-\sqrt{3} \\ -\sqrt{3}
	\end{bmatrix}. 
\end{aligned}
\]
\second{We use the same approach presented in \cite{calzola23} to treat the homogeneous Neumann boundary conditions}. For $y \notin \Omega_1\times\Omega_2$ we define $\pi^{\vec{n}}(\cdot)$, the projection operator along the direction $\vec{n}$ orthogonal to the boundary $\partial\left(\Omega_1\times\Omega_2\right)$. Whenever a characteristic curve $y_{i,j}^{n,\ell}$, for some values of $i,j,n$ and $\ell$, falls out of the domain, we reflect it inside in the point
\[
\pi^{\vec{n}}(y_{i,j}^{n,\ell}) - \overline c \sqrt{\Delta t}\vec{n}(\pi^{\vec{n}}(y_{i,j}^{n,\ell}))
\]
where $\overline c > 0$ is a fixed constant \second{that should be chosen is such a way that the reflected characteristic stays inside the domain and the distance between the original characteristic and the reflected one is of order $\sqrt{\Delta t}$. This last request is the fundamental assumption to obtain stability of the scheme (for more information, we refer the reader to \cite[Remarks 3 and 5]{calzola23})}. This leads to the modified definition of the characteristic curves
\begin{equation}\label{eq:reflchar}
	\tilde{y}_{i,j}^{n,\ell} = \begin{cases} y_{i,j}^{n,\ell}, \quad&\mbox{if } y_{i,j}^{n,\ell}\in\Omega_1\times\Omega_2,\\
		\pi^{\vec{n}}(y_{i,j}^{n,\ell}) - \overline c \sqrt{\Delta t}\vec{n}(\pi^{\vec{n}}(y_{i,j}^{n,\ell})), &\mbox{if } y_{i,j}^{n,\ell}\notin\Omega_1\times\Omega_2.
	\end{cases}
\end{equation}
We introduce the following notation for the reaction term in \eqref{eq:hj}
\begin{equation*}%\label{eq:react}
	\begin{aligned}
		\mathcal{R}[f,\psi,\mathbf{u}](\mathbf{v},t) = \mathcal{Q}_1[f,\psi](\mathbf{v},t) + \mathcal{Q}_2[f,\psi](\mathbf{v},t) \\
		+ \frac{1}{2}\left(\lvert v_1-\overline v_1\rvert^2s_1(v) + \lvert v_2-\overline v_2\rvert^2s_2(v) + \gamma \lVert \mathbf{u}\rVert^2 \right),
	\end{aligned}
\end{equation*}
so that the scheme for the backward equation reads as
\begin{equation}\label{eq:sl_hj}
	\begin{aligned}
		\psi_{i,j}^n = &\sum_{\ell=1}^9\omega^\ell\left(\mathbb{I}^r[\pmb{\psi}^{n+1}]\left(\tilde{y}_{i,j}^{n,\ell}\right) + \frac{\Delta t}{2} \mathbb{I}^r\left[\mathcal{R}[\mathbf{f}^{n+1},\pmb{\psi}^{n+1},\mathbf{u}^{n+1}]\right](\tilde{y}_{i,j}^{n,\ell})\right) \\
		&+ \frac{\Delta t}{2} \mathcal{R}[\mathbf{f}^{n},\pmb{\psi}^{n},\mathbf{u}^{n}](\mathbf{v}_{i,j},t^{n}),
	\end{aligned}
\end{equation}
for all $n=0,\dots,N_T-1$ and with $\omega^1 = 4/9$, $\omega^2=\omega^3=\omega^4=\omega^5=1/9$, $\omega^6=\omega^7=\omega^8=\omega^9=1/36$. Since we are dealing with homogeneous Neumann boundary conditions, the correction term appearing in the scheme described in \cite{calzola23} is equal to zero. \second{The only modification with respect to a classical second-order semi-Lagrangian discretization consists in the use of the reflected characteristic defined in \eqref{eq:reflchar}}.
\subsection{The \second{structure-preserving} scheme for the forward Fokker--Planck equation}\label{sec:spscheme}
For the forward equation in \eqref{eq:opifp}\second{,} we apply the \second{structure-preserving} scheme in \cite{pareschi18}, which is in practice a Chang--Cooper type of scheme. To take advantage of the \second{absence of a parabolic CFL condition for the stability of the semi-Lagrangian method}, we implemented a semi-implicit version of such scheme. We choose the same discretization steps $\Delta v_1$,$\Delta v_2$ and $\Delta t$ that we used for the semi-Lagrangian scheme. We rewrite equation \eqref{eq:opifp} in divergence form
\begin{equation}
	\partial_t f = \nabla \cdot \mathcal{F}[f,u](v_1,v_2,t),
\end{equation}
with
\begin{equation}\label{eq:flux}
	\mathcal{F}[f,u](v_1,v_2,t) =\left(-\mathcal{G}[f;u](v_i,v_2,t) + \nabla \cdot \mathcal{D}(v_1,v_2)\right)f + \mathcal{D}(v_1,v_2)\nabla f
\end{equation}
Using the uniform spatial grid previously introduced, we denote \second{by} \first{$\mathbf{v}_{i\pm\frac{1}{2},j} = \mathbf{v}_{i,j}\pm\Delta v_1/2\mathbf{\hat v}_1$ and $\mathbf{v}_{i,j\pm\frac{1}{2}} = \mathbf{v}_{i,j}\second{\pm}\Delta v_2/2\mathbf{\hat v}_2$, with 
	\[\mathbf{\hat v}_1 = \begin{bmatrix}1 \\ 0\end{bmatrix}, \quad \mbox{and} \quad \mathbf{\hat v}_2 = \begin{bmatrix} 0 \\ 1 \end{bmatrix}
	\]
	unit vectors of the two axis}. Supposing the control $\mathbf{u}$ to be given, we use the notation $\mathcal{F}_{i,j}^{(1)}$ and $\mathcal{F}_{i,j}^{(2)}$, respectively, for the first and the second component of the flux \eqref{eq:flux} in $\mathbf{v}_{i,j}$. We consider the semi-discrete approximation
\begin{equation}\label{eq:fp_semid}
	\frac{\mathrm{d}}{\mathrm{d}t}f_{i,j}(t) = \frac{\mathcal{F}^{(1)}_{i+\frac{1}{2},j}(t) - \mathcal{F}^{(1)}_{i-\frac{1}{2},j}(t)}{\Delta v_1} + \frac{\mathcal{F}^{(2)}_{i,j+\frac{1}{2}}(t) - \mathcal{F}^{(2)}_{i,j-\frac{1}{2}}(t)}{\Delta v_2},
\end{equation}
for all $t\geq 0$, where $\mathbf{f} = (f_{i,j})_{i,j}$, $f_{i,j}(t)$ approximates $f(\first{\mathbf{v}_{i,j}},t)$, $\mathcal{F}^{(1)}_{i\pm\frac{1}{2},j}$ is the numerical flux along the first axis \second{and} $\mathcal{F}^{(2)}_{i,j\pm\frac{1}{2}}$ is the numerical flux along the second axis. Defining $\mathcal{A}[f;\mathbf{u}](\first{\mathbf{v}},t) = -\mathcal{P}_1[f](\first{\mathbf{v}},t) - u_1(\first{\mathbf{v}},t) + \partial_{v_1}\mathcal{D}_{1}(\first{\mathbf{v}})$ and $\mathcal{B}[f;\mathbf{u}](\first{\mathbf{v}},t) = -\mathcal{P}_2[f](\first{\mathbf{v}},t) - u_2(\first{\mathbf{v}},t) + \partial_{v_2}\mathcal{D}_{2}(\first{\mathbf{v}})$, we consider the fluxes
\begin{equation}\label{eq:fluxes}
	\begin{aligned}
		\mathcal{F}^{(1)}_{i+\frac{1}{2},j}\first{[f]} &= \first{\mathcal{A}_{i+\frac{1}{2},j}[f]}\tilde f_{i+\frac{1}{2},j} + \mathcal{D}_{1}(\first{\mathbf{v}_{i+\frac{1}{2},j}})\frac{f_{i+1,j}-f_{i,j}}{\Delta v_1},\\
		\mathcal{F}^{(2)}_{i,j+\frac{1}{2}}\first{[f]} &= \first{\mathcal{B}_{i,j+\frac{1}{2}}[f]}\tilde f_{i,j+\frac{1}{2}} + \mathcal{D}_{2}(\first{\mathbf{v}_{i,j+\frac{1}{2}}})\frac{f_{i,j+1} - f_{i,j}}{\Delta v_2},
	\end{aligned}
\end{equation}
in which we used the notations
\[
\tilde f_{i+\frac{1}{2},j} = \left(1-\delta_{i+\frac{1}{2},j}\right)f_{i+1,j} +\delta_{i+\frac{1}{2},j}f_{i,j},
\]
and
\[
\tilde f_{i,j+\frac{1}{2}} = \left(1-\delta_{i,j+\frac{1}{2}}\right)f_{i,j+1} +\delta_{i,j+\frac{1}{2}}f_{i,j},
\]
\first{and set $\mathcal{A}_{i+\frac{1}{2},j}[f] = \mathcal{A}[f;\mathbf{u}](\first{\mathbf{v}_{i+\frac{1}{2},j}},t)$ and $\mathcal{B}_{i,j+\frac{1}{2}}[f] = \mathcal{B}[f;\mathbf{u}](\first{\mathbf{v}_{i,j+\frac{1}{2}}},t)$.}
The choice of the value of the coefficients $\delta_{i+\frac{1}{2},j}$ and $\delta_{i,j+\frac{1}{2}}$ is made in such a way that the stationary state structure is preserved (we refer to \cite{pareschi18} for a detailed explanation on how this can be obtained). In particular,
\[
\delta_{i+\frac{1}{2},j} = \frac{1}{\lambda_{i+\frac{1}{2},j}}+\frac{1}{1-\exp\left\{\lambda_{i+\frac{1}{2},j}\right\}},\]
\[\mbox{with}\quad \lambda_{i+\frac{1}{2},j} = \frac{\Delta v_1\left(\mathcal{A}_{i+\frac{1}{2},j}\first{[f]}+\partial_{v_1}\mathcal{D}_{1}(\first{\mathbf{v}_{i+\frac{1}{2},j}})\right)}{\mathcal{D}_{1}\left(\first{\mathbf{v}_{i+\frac{1}{2},j}}\right)},
\]
and
\[
\delta_{i,j+\frac{1}{2}} = \frac{1}{\lambda_{i,j+\frac{1}{2}}}+\frac{1}{1-\exp\left\{\lambda_{i,j+\frac{1}{2}}\right\}},\]
\[\mbox{with}\quad \lambda_{i,j+\frac{1}{2}} = \frac{\Delta v_2\left(\mathcal{B}_{i,j+\frac{1}{2}}\first{[f]}+\partial_{v_2}\mathcal{D}_{2}(\first{\mathbf{v}_{i,j+\frac{1}{2}}})\right)}{\mathcal{D}_{2}(\first{\mathbf{v}_{i,j+\frac{1}{2}}})},
\]
a choice that guarantees second-order accuracy in space. We now introduce the same time discretization used for the semi-Lagrangian scheme and use a second-order in time discretization for the temporal derivatives appearing in \eqref{eq:fp_semid}. In particular, we use a second-order implicit-explicit (IMEX) Runge-Kutta scheme defined as follows (\cite{pareschi18,boscarino16}). \first{The choice of relying on an IMEX scheme is dictated here by the need to avoid a parabolic CFL for the positivity of this scheme, as proven in \cite[Propositions 2 and 3]{pareschi18}. A fully explicit implementation would require $\Delta t = O(\Delta x)$, a restriction that would invalidate the main advantage of the semi-Lagrangian scheme for the backward equation.} We formally duplicate the number of variables so that our problem becomes
\begin{equation}
	\begin{aligned}
		\frac{\mathrm{d}}{\mathrm{d}t}f_{i,j}(t) &=\,\mathcal{S}_{i,j}[\mathbf{f},\mathbf{g}], \\
		\frac{\mathrm{d}}{\mathrm{d}t}g_{i,j}(t) &=\, \mathcal{S}_{i,j}[\mathbf{f},\mathbf{g}],
	\end{aligned}
\end{equation}
where
\[
\mathcal{S}_{i,j}[\mathbf{f},\mathbf{g}] = \frac{\mathcal{F}^{(1)}_{i+\frac{1}{2},j}[\mathbf{f},\mathbf{g}] (t) - \mathcal{F}^{(1)}_{i-\frac{1}{2},j}[\mathbf{f},\mathbf{g}](t)}{\Delta v_1} + \frac{\mathcal{F}^{(2)}_{i,j+\frac{1}{2}}[\mathbf{f},\mathbf{g}](t) - \mathcal{F}^{(2)}_{i,j-\frac{1}{2}}[\mathbf{f},\mathbf{g}](t)}{\Delta v_2},
\]
and
\begin{equation}\label{eq:fluxes_rk}
	\begin{aligned}
		\mathcal{F}^{(1)}_{i+\frac{1}{2},j}[\mathbf{f},\mathbf{g}] =& \mathcal{A}_{i+\frac{1}{2},j}[\mathbf{g}]\left(\left(1-\delta_{i+\frac{1}{2},j}[\mathbf{g}]\right)f_{i+1,j} +\delta_{i+\frac{1}{2},j}[\mathbf{g}]f_{i,j}\right) \\ 
		&+ \first{\mathcal{D}}_{1}(v_{i+\frac{1}{2},j})\frac{f_{i+1,j}-f_{i,j}}{\Delta v_1},\\
		\mathcal{F}^{(2)}_{i,j+\frac{1}{2}}[\mathbf{f},\mathbf{g}] =& \mathcal{B}_{i,j+\frac{1}{2}}[\mathbf{g}]\left(\left(1-\delta_{i,j+\frac{1}{2}}[\mathbf{g}]\right)f_{i,j+1} +\delta_{i,j+\frac{1}{2}}[\mathbf{g}]f_{i,j}\right) \\
		&+ \first{\mathcal{D}}_{2}(v_{i,j+\frac{1}{2}})\frac{f_{i,j+1} - f_{i,j}}{\Delta v_2}.
	\end{aligned}
\end{equation}
where $[\cdot]$ denotes functional dependence. Defining $g_{i,j}^n = f_{i,j}^n$, the implicit-explicit Runge-Kutta scheme can be written as
\begin{equation}
	\begin{cases}
		\displaystyle f_{i,j}^{n+1} = f_{i,j}^n + \Delta t \sum_{k = 1}^s b_k \mathcal{S}_{i,j}[\mathbf{F}^k,\mathbf{G}^k], \\[6pt]
		\displaystyle g_{i,j}^{n+1} = f_{i,j}^n + \Delta t \sum_{k = 1}^s \tilde{b}_k \mathcal{S}_{i,j}[\mathbf{F}^k,\mathbf{G}^k],
	\end{cases}
\end{equation}
where $\mathbf{F}^h$ and $\mathbf{G}^h$, $h = 1, \dots,s$ are defined by
\begin{equation}
	\begin{cases}
		\displaystyle F_{i,j}^{h} = f_{i,j}^n + \Delta t \sum_{k = 1}^h a_{hk} \mathcal{S}_{i,j}[\mathbf{F}^k,\mathbf{G}^k], \\[6pt]
		\displaystyle G_{i,j}^{h} = f_{i,j}^n + \Delta t \sum_{k = 1}^{h-1} \tilde{a}_{hk} \mathcal{S}_{i,j}[\mathbf{F}^k,\mathbf{G}^k],
	\end{cases}
\end{equation}
for all $h = 1, \dots, s$. In order to make the duplication of the system only apparent we set $b_k = \tilde{b}_k$, so $f_{i,j}^{n+1} = g_{i,j}^{n+1}$, and to obtain second-order accuracy in time we set $s=2$ and
\[
a_{11} = 0,\quad a_{21} = a_{22} = \frac{1}{2},\quad \tilde{a}_{21} = 1,\quad \mbox{and}\quad b_k = \tilde{b}_k = \frac{1}{2}, \quad k = 1,2.
\]
\begin{remark}
	\first{Using a Chang-Cooper type of scheme lets us naturally treat the zero-flux boundary conditions associated with the Fokker--Planck equation. In fact, it is sufficient to set the flux at the boundary equal to zero, i.e.
		\[
		\mathcal{F}_{-\frac{1}{2},j} = \mathcal{F}_{N_{v_1}+\frac{1}{2},j}\quad \mbox{and} \quad \mathcal{F}_{i,-\frac{1}{2}} = \mathcal{F}_{i,N_{v_2}+\frac{1}{2}} = 0,
		\]
		for all $i = 0, \dots, N_{v_1}-1$ and $j = 0,\dots,N_{v_2}-1$, to guarantee that no part of the boundary is outflow.}
	
	\first{Furthermore, we stress that the Chang-Cooper scheme is particularly designed to capture the stationary states of Fokker-Planck equation. Indeed, considering the 1D setting for the sake of simplicity, once we imposed the steady state on each interface, we have the following relation
		\[
		\frac{f_{i+1}}{f_i} = \exp\{\lambda_{i+1}\},
		\]
		where the coefficients $\lambda_{i+1/2}$ are designed as follows
		\[
		\lambda_{i+1/2} := \int_{v_i}^{v_{i+1}} \frac{\partial_v \mathcal{D}(v)-\mathcal{P}[f](v,t)-u(v,t)}{\mathcal{D}(v)}\,dv.
		\]
		Thus, it is evident that only for simplified cases, such as linear drifts with constant control input, stationary states can be captured exactly, by exact integration. Otherwise, for nonlinear drifts and in presence of time varying control steady states can be only accurately approximated on each cell, where the accuracy depends on the quadrature rule employed for the weights $\lambda_{i+1/2}.$ We refer to \cite{pareschi18,chang1970practical} for further details on the properties of such schemes for Fokker-Planck equations.} 
\end{remark}
\second{We provide a schematic algorithmic representation of our method in order to clarify the main steps.}

\begin{algorithm}
	\caption{Numerical solution to optimal control problem \eqref{eq:minimization}-\eqref{eq:opifp}.}\label{alg:cap}
	
	\begin{algorithmic}
		\State $f^{(0)} \gets f_0$ for all $n$\vspace{2pt}
		\State $\psi^{(0)} \gets 0$ for all $n$\vspace{2pt}
		\State $\mathbf{u}^{(0)} \gets 0$ for all $n$\vspace{2pt}
		\State $\lambda^{(0)} \gets \overline \lambda$\vspace{2pt}
		\State $\verb|err_cost| \gets 1$\vspace{2pt}
		\State $N \gets 1$\vspace{2pt}
		\While{$\texttt{err{\_}cost} > 10^{-5}$ and $N < 500$}\vspace{2pt}
		\State $\psi^{(N)} \gets \verb|SL_scheme|(f^{(N-1)},\psi^{(N-1)},\mathbf{u}^{(N-1)})$\vspace{2pt}
		\State $f^{(N)} \gets \verb|ChangCooper_scheme|(\psi^{(N)},\mathbf{u}^{(N-1)})$\vspace{2pt}
		\If{$N>1$}\vspace{2pt}
		\State $\texttt{DL} \gets \left(D_u\mathcal{L}(f^{(N)},\mathbf{u}^{(N-1)},\psi^{(N)}) - D_u\mathcal{L}(f^{(N-1)},\mathbf{u}^{(N-2)},\psi^{(N-1)})\right)$\vspace{2pt}
		\State $\displaystyle\lambda^{(N)} \gets \frac{\lvert \left(\mathbf{u}^{(N-1)} - \mathbf{u}^{(N-2)}\right) \cdot \texttt{DL} \rvert}{\lVert \texttt{DL} \rVert^2}$\vspace{2pt}
		\EndIf\vspace{2pt}
		\State $\mathbf{u}^{(N)} \gets \mathbf{u}^{(N-1)} - \lambda^{(N)}D_u\mathcal{L}(f^{(N)},\mathbf{u}^{(N-1)},\psi^{(N)})$\vspace{2pt}
		\State $\texttt{err{\_}cost} \gets \lvert \mathcal{J}(\mathbf{u}^{(N)};f_0) - \mathcal{J}(\mathbf{u}^{(N-1)};f_0)\rvert$\vspace{2pt}
		\State $N \gets N+1$\vspace{2pt}
		\EndWhile
	\end{algorithmic}
\end{algorithm}

\section{Numerical simulations}\label{sec:numerics}
In \second{this Section we report some numerical simulations to test the accuracy of the proposed methodology}. Section \ref{sec:num1} and \ref{sec:num2} \second{are dedicated to showing the accuracy order of the schemes: in the simulation provided,} given $h$ the solution obtained with a numerical scheme, we compute the normalized $L^2$ and $L^\infty$ errors with respect to an exact (or reference) solution $\hat h$, defined as
\begin{equation}\label{eq:err}
	E_2(h) = \frac{\lVert h-\hat h \rVert_2}{\lVert \hat h \rVert_2}, \quad E_\infty(h) = \frac{\lVert h-\hat h \rVert_\infty}{\lVert \hat h \rVert_\infty}.
\end{equation}
\second{Sections \ref{sec:opic} and \ref{sec:dueopi} present instead two tests in two spatial dimensions: the first one can be understood as a problem of opinion control when a connectivity structure is present, while the second one can be seen as a problem in which two interacting opinions of a population are controlled. All numerical simulations were performed on a laptop equipped with an Intel(R) Core(TM) i7-10750H CPU @ 2.60GHz, 16 GB of RAM, and a 64-bit operating system on an x64-based processor. A minimal reference code as the one in Section \ref{sec:compo1o2} is available on GitHub: \href{https://github.com/ClzLse/order2controlFP}{ClzLse/order2controlFP}.}
\second{\begin{remark}
		We must stress the fact that, even though both the semi-Lagrangian and the Chang-Cooper type scheme do not require a parabolic CFL condition in order to be stable, one cannot simply choose any ratio between the time and the space discretization steps. In fact, in order to have a smooth domain of dependence, the semi-Lagrangian scheme in presence of diffusion must verify a compatibility condition
		\[
		\frac{\max(\mathcal{D})^{1/3}\Delta t}{T^{2/3}\Delta x^{2/3}}\ll 1
		\]
		(see \cite{ferretti10}), while the ratio between $\Delta t$ and $\Delta x$ for the structure-preserving scheme depends on the maximum value of the drift, that in this case contains the control itself. This means that, in presence of very high diffusion, the method becomes slower because the compatibility condition of the semi-Lagrangian scheme implies a limitation on the value of $\Delta t$. Instead, if the drift term, especially due to the control, becomes large, then the limitation is due to the condition for the positivity of the solution of the structure-preserving scheme. 
	\end{remark}
}
\subsection{Convergence to the stationary state for the Fokker--Planck}\label{sec:num1}
In this first simulation\second{,} we show the second-order convergence of the solution to the Fokker--Planck to a stable stationary state. The general expression of the steady state of \eqref{eq:opifp} is not known, \second{and multiple stationary states may exist due to the presence of the non-local operator $\mathcal{P}[f](\cdot)$ and the nonlinear diffusion coefficient $\mathcal{D}(\cdot)$, see for example \cite{frank2005nonlinear}.}

\begin{figure}[h!]
	\centering
	\includegraphics[width=0.45\textwidth]{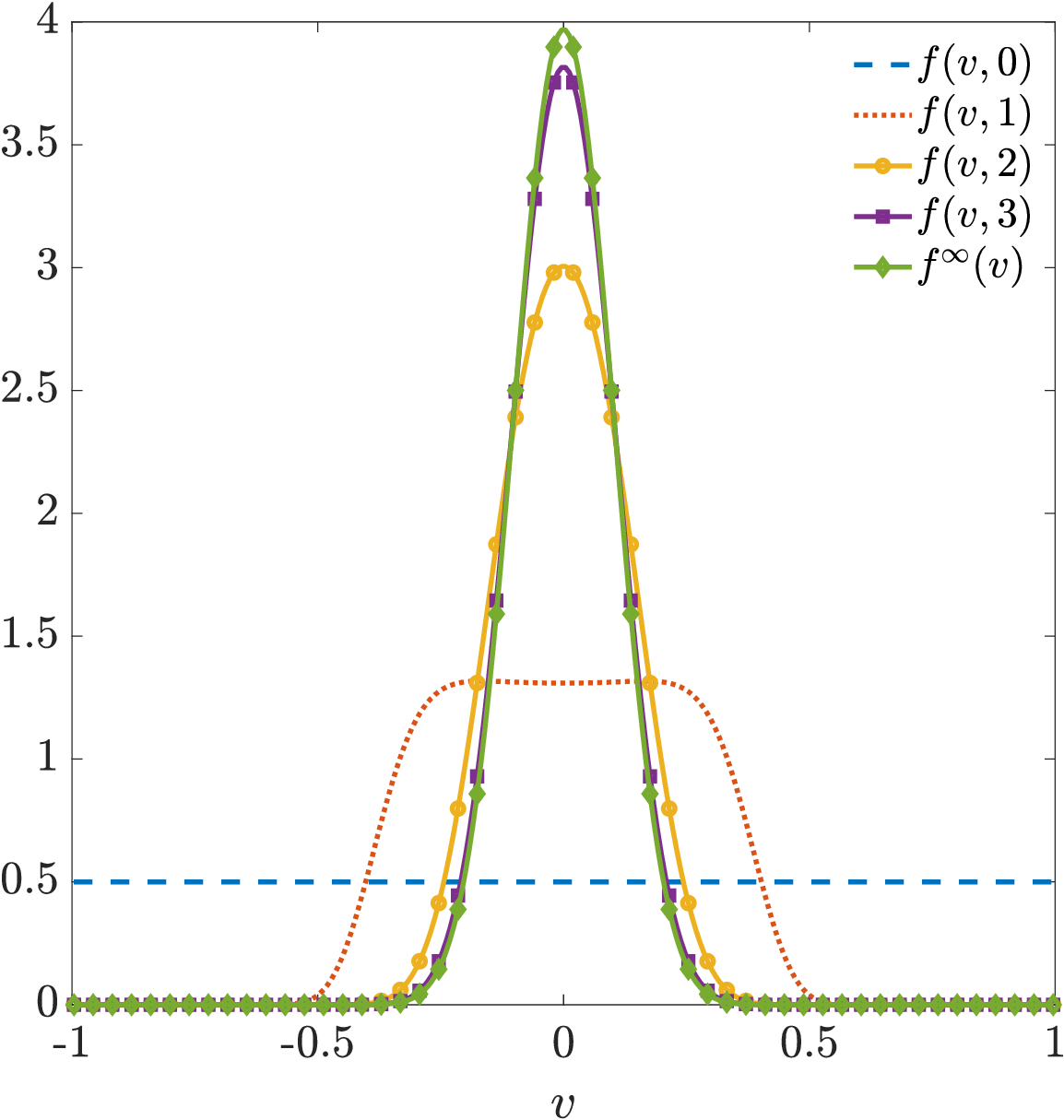}\quad\includegraphics[width=0.46\textwidth]{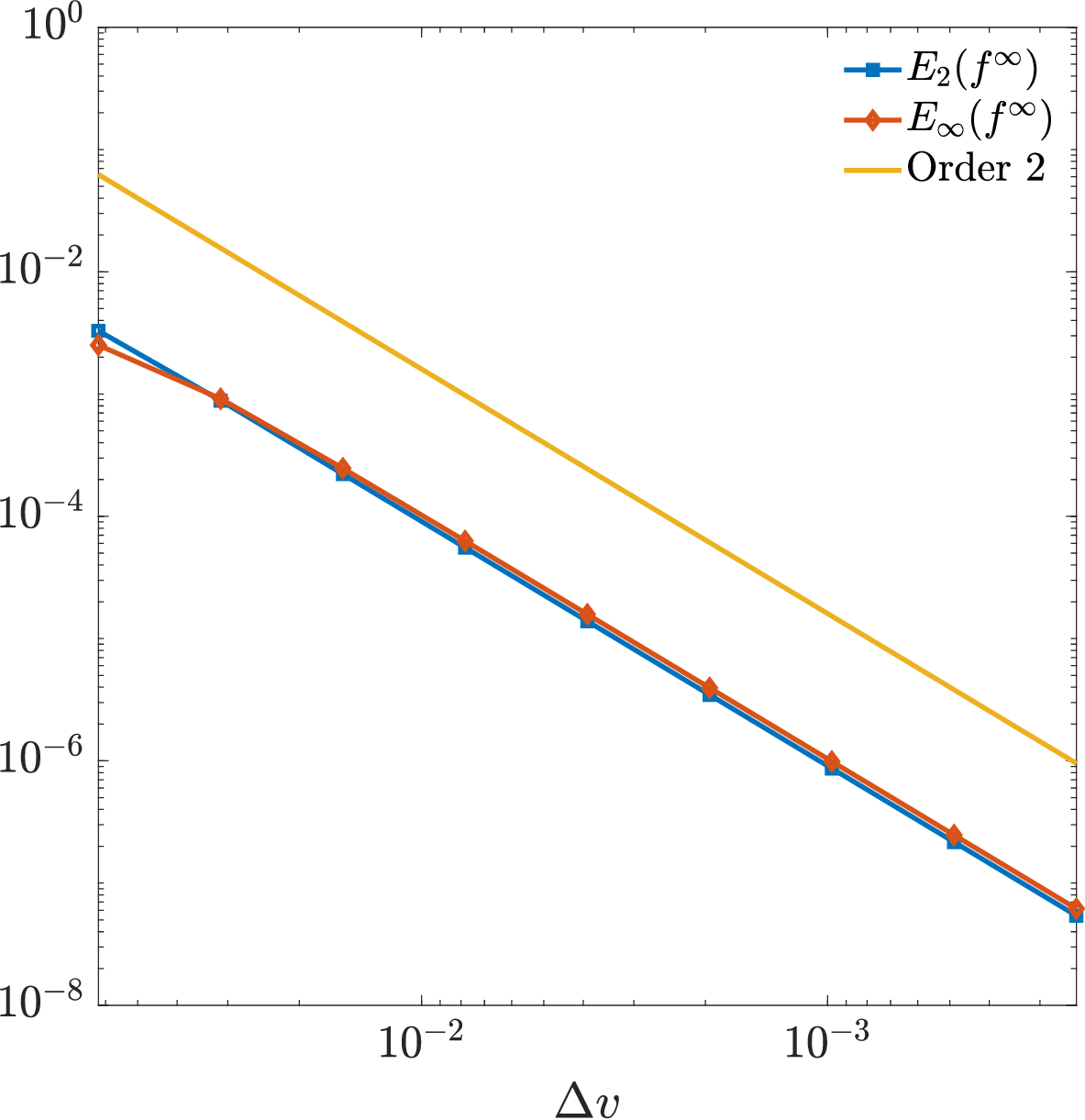}
	\caption{Time evolution of the numerical approximation of $f(v,t)$ (left) and convergence to the stationary state (right). The direction of the horizontal axis is inverted in the second picture.}
	\label{fig:test_staz}
\end{figure}

Nevertheless, under some simplifying assumptions it is possible to compute such stationary states. We put ourselves in the one dimensional setting, $v\in\Omega=[-1,1]$, and we suppose that $\first{\mathcal{D}}(v) = \frac{\sigma^2}{2}(1 - v^2)^2$, $u(v,t) = 0$, and $P(v,v_*)=1$ so that the nonlocal drift functional becomes
\begin{equation*}
	\mathcal{P}[f](v,t) = \displaystyle\left(\int_{\Omega} v_* f(v_*,t)\,dv_* - v \right) 
	%\\
	%&=& \alpha \left(\int_{I\times\R_+} v_* g(v_*, \tau)\,dv_* - v \right)
	=\left( m_v(t) - v \right),
\end{equation*}
\second{where we used the notation $m_v(t)$ to indicate the mean of $f$ at time $t$, i.e. $m_v(t) = \int_\Omega v_*f(v_*,t)\mathrm{d}v_*$.}
%We also assume, in this situation, that we can look for solutions of type $f(v,c,t) = f_1(v,t)f_2(c,t)$ that leads to asymptotic steady states of the form $f^\infty(v,c)=f_{1,\infty}(v) f_{2,\infty}(c)$. 
The  Fokker--Planck equation \eqref{eq:opifp} can then be rewritten as

\begin{equation}\label{eq:fpeq2}
	\frac{\partial f}{\partial t} = - \frac{\partial ( (m_v(t)  - v)f)}{\partial v}+ \frac{\sigma^2}{2} \frac{\partial^2 ((1-v^2)^2f)}{\partial v^2},
\end{equation}
whose 
%solutions can be found by solving
%\begin{equation}\label{eq:perg}
%	\frac{\partial f_1}{\partial t} = - \alpha\frac{\partial ( (m_v(t) - v) f_1)}{\partial v} + \frac 12 \sigma^2 \frac{\partial^2 ((1-v^2)^2f_1)}{\partial v^2} 
%\end{equation} 
%and
%\begin{equation}\label{eq:perh}
%	\frac{\partial f_2}{\partial t} = \frac{\mu}{2}  \frac{\partial \left(c \text{ln} \left(\frac{c}{\bar c}\right)f_2\right)}{\partial c} + \frac 12 \nu^2  \frac{\partial^2( c^2f_2)}{\partial c^2}.
%\end{equation} 
%The solution to \eqref{eq:perh} is a log-normal distribution, while the 
stationary solution $f^\infty(v)$ is given by
$$
f^{\infty}(v) = K_\infty(1+v)^{-2+\bar m_v/2\sigma^2}(1-v)^{-2-\bar m_v/2\sigma^2}\text{exp}\left\{ -\frac{(1-\overline m_v v)}{\sigma^2(1-v^2)}\right\},
$$
where \second{$\lim_{t\to \infty} m_v(t) = \overline m_v = m_v(0)$ since the mean is conserved, and} $K_\infty>0$ is such that 
\[ \int_{\Omega} f^{\infty}(v) \mathrm{d}v = 1 \]
(we refer, for instance, to \cite{toscani06,albi24} for more details on the computations). We focus on the convergence of the numerical solution to \eqref{eq:fpeq2} using the scheme described in Section \ref{sec:spscheme}. Starting from an initial condition given by the uniform distribution in $\Omega$, we check the rate of convergence of our scheme to the analytical steady state with $\sigma^2 = 2\cdot 10^{-2}$, and $\Delta t = \Delta v/2$. \first{We consider the steady state to be reached when $\max\lvert f^{n+1}-f^n\rvert < \texttt{toll}$, with $\texttt{toll} = 10^{-12}$}. \second{We define $\Delta v = \lvert \Omega \rvert/2^{N_v}$, $N_v\in\mathbb{N}$, and  we compute the errors for several choices of $N_v$}. Figure \ref{fig:test_staz}, on the left, shows the initial condition, the numerical solution at time $t=1,2,3$ and the stationary state for $N_v = 12$, while on the right we have the \second{errors $E_2$ and $E_\infty$} defined in \eqref{eq:err} compared to \second{a straight line with slope equal to $2$}. Table \ref{tab_staz} reports the value of the errors for each choice of numbers of the intervals $2^{N_v}$ used in the discretization of $\Omega$. \second{The values of $p_2$ and $p_\infty$ indicate the order of convergence with respect to errors $E_2$ and $E_\infty$, respectively. The expected value of $2$ is reached for both norms.}

\begin{table}[ht]
	\centering\begin{minipage}{.8\textwidth}
		\caption{Errors and convergence rates for the Fokker--Planck equation at the stationary state. In all simulations $\Delta t = \Delta v/2$.}
		\label{tab_staz}
		\begin{tabular*}{\textwidth}{@{\extracolsep{\fill}}*{6}{l}@{\extracolsep{\fill}}}
			\toprule%
			$N_v$&$E_2(f^{\infty})$&$p_2$&$E_\infty(f^{\infty})$&$p_\infty$&\second{CPU time $(s)$}\\
			\midrule
			$4$&$3.2925\times 10^{-3}$&$/$&$2.5138\times 10^{-3}$&$/$&$3.30\times 10^{-2}$ \\
			\midrule
			$5$&$8.8520\times 10^{-4}$&$1.90$&$9.1545\times 10^{-4}$&$1.46$&$5.82\times 10^{-2}$ \\
			\midrule
			$6$&$2.2177\times 10^{-4}$&$2.00$&$.4692\times 10^{-4}$&$1.89$&$1.27\times 10^{-1}$\\
			\midrule
			$7$&$5.5473\times 10^{-5}$&$2.00$&$6.2862\times 10^{-5}$&$1.97$&$4.79\times 10^{-1}$ \\
			\midrule
			$8$&$1.3870\times 10^{-5}$&$2.00$&$1.5786\times 10^{-5}$&$1.99$& $2.10$\\
			\midrule
			$9$&$3.4676\times 10^{-6}$&$2.00$&$3.9510\times 10^{-6}$&$2.00$&$2.04\times 10^{1}$\\
			\midrule
			$10$&$8.6691\times 10^{-7}$&$2.00$&$9.8802\times 10^{-7}$&$2.00$& $1.22\times 10^{2}$\\
			\midrule
			$11$&$2.1673\times 10^{-7}$&$2.00$&$2.4702\times 10^{-7}$&$2.00$&$1.09\times 10^{3}$\\
			\midrule
			$12$&$5.4182\times 10^{-8}$&$2.00$&$6.1757\times 10^{-8}$&$2.00$&$8.92\times 10^{3} $\\
			\bottomrule
	\end{tabular*}\end{minipage}
\end{table}

\subsection{Order check with fixed suboptimal control for the Hamilton--Jacobi equation}\label{sec:num2}
\begin{figure}[h!]
	\centering
	\includegraphics[width=0.45\textwidth]{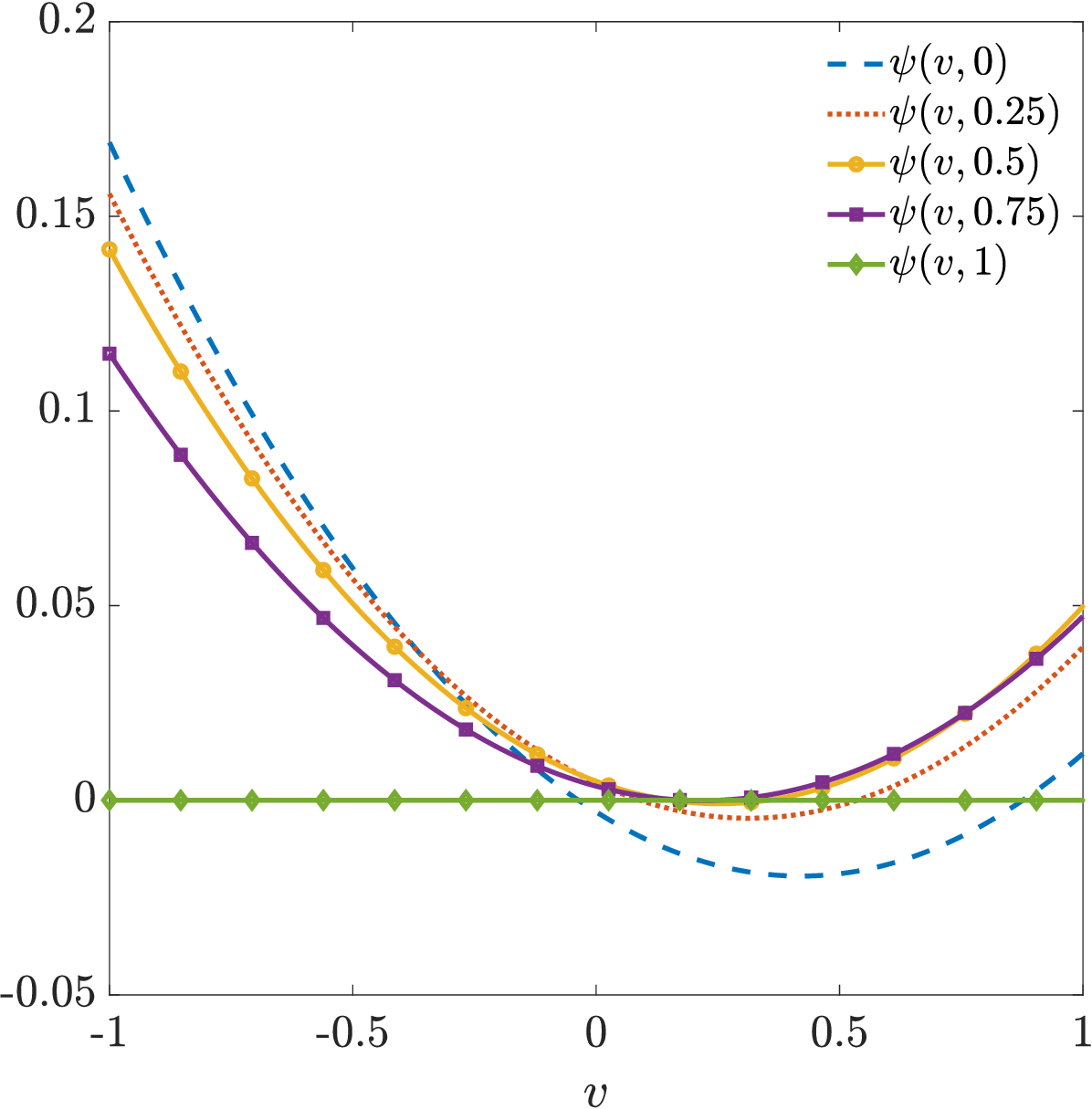}\quad\includegraphics[width=0.445\textwidth]{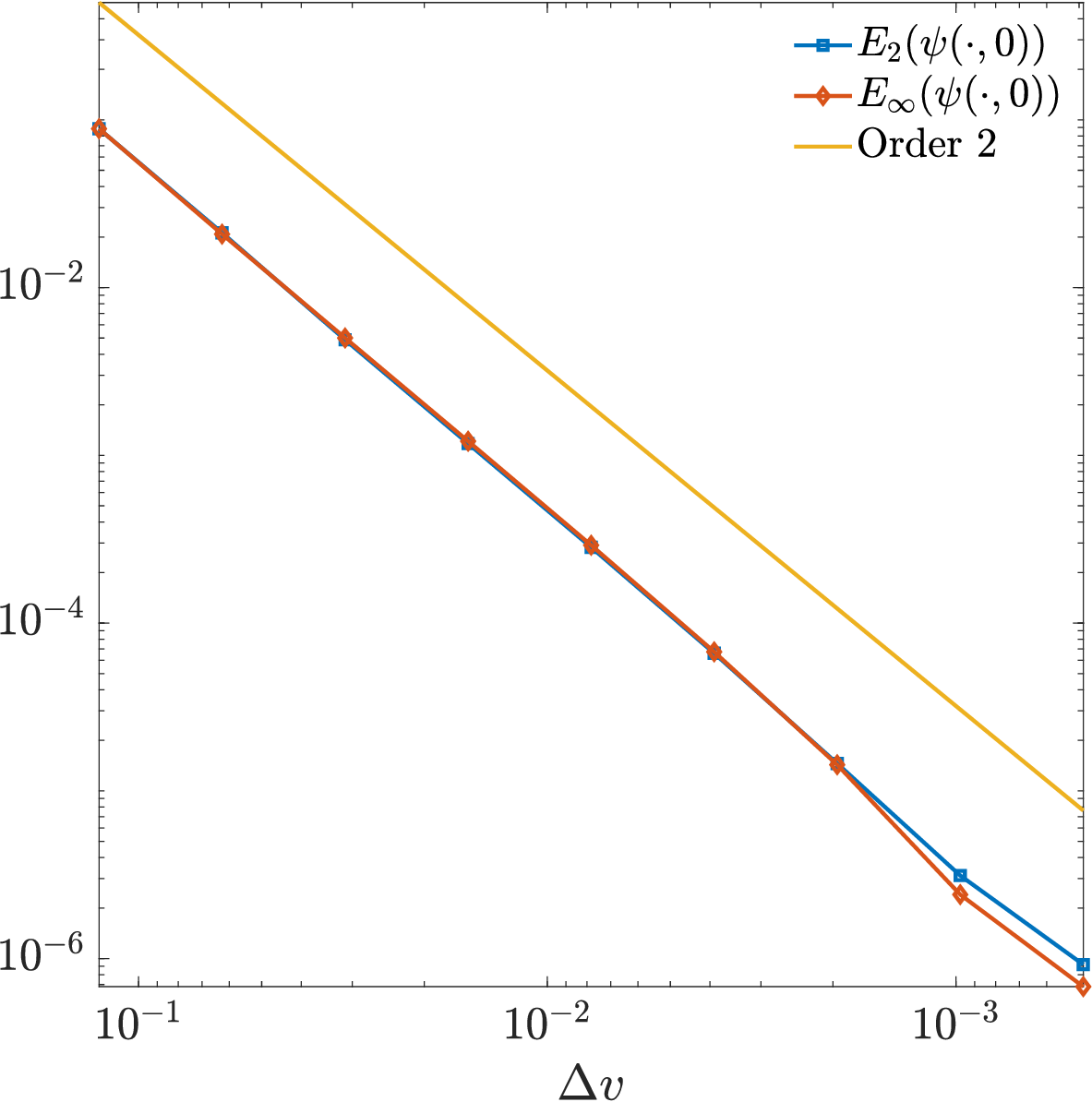}
	\caption{Time evolution of the numerical approximation of $\psi(v,t)$ (left) and convergence to the reference solution (right). The direction of the horizontal axis is inverted in the second picture.}
	\label{fig:test_staz_psi}
\end{figure}

Consider again the one dimensional setting. We want to simulate the solution to \eqref{eq:hj} with a fixed $f$ and a suboptimal control, in order to check the second-order accuracy of the semi-Lagrangian scheme with nonlocal drift and forcing term. First, we compute a reference solution on a spatial grid with $N_v = 2^{13}$ discretization points solving 
\begin{equation}
	\begin{cases}
		\begin{aligned}
			\partial_t f(v,t) &=-\first{\nabla \cdot}\left(\left(\mathcal{P}[f](v,t) + u(v,t)\right) f -\first{\nabla}\left(\first{\mathcal{D}}(v)f\right)\right), \,\, &(v,t) \in \Omega \times (0,T]\\
			f(v,0) &= K_0\exp{(-(v-\lambda)^2/(2\rho^2))}, &v\in \Omega\\
			
			-\partial_t \psi &=\left(\mathcal{P}[f](v,t) + u(v,t)\right)\partial_v \psi + \first{\mathcal{D}}(v)\first{\nabla^2}_{v} \psi &\\ 
			& \quad+ \mathcal{Q}[f,\psi](v,t) + \frac{1}{2}\left(\vert v-\overline v \vert^2 s(c) + \gamma \vert u \vert^2 \right), &(v,t) \in \Omega \times [0,T),\\
			\psi(v,\first{T}) &= 0, &v\in \Omega,
		\end{aligned}
	\end{cases}
\end{equation}
coupled with zero-flux boundary conditions for the Fokker--Planck equation and homogeneous Neumann boundary conditions for the Hamilton--Jacobi. We set $P(v,v_*) = 1$, $h(v)=0$ in \eqref{eq:drift_general}, and \second{the} diffusion equal to $\first{\mathcal{D}}(v) = \sigma^2(1-v^2)^2/2$ with  $\sigma^2 = 1\times 10^{-2}$. \second{The final time is set to $T=1$. In \eqref{eq:hj} we \second{set} $s(c)=1$, $\overline v=0.2$, $\gamma = 0.05$}, and we fix as suboptimal control 
\[
u(v,t) =\left(-\frac{5}{2}v + \frac{1}{2}\right)(T-t).
\]
The initial condition for the Fokker--Planck is given by 
\[
f_0(v) = K_0e^{-\frac{(v+0.5)^2}{0.5}}
\]
\second{i.e. $f_0$ is a normal distribution with mean equal to $-0.5$ and variance $0.25$}, where the constant $K_0>0$ is chosen so that 
\[
\int_\Omega f_0(v)\mathrm{d}v = 1.
\]
Once we have the reference solution for both the Fokker--Planck and the Hamilton--Jacobi\second{,} we use the reference $f(\first{v},t)$ as a datum and compute the errors for equation \eqref{eq:hj} at time $t=0$. Table \ref{tab2} report\second{s} the errors as defined in \eqref{eq:err} and confirms the order $2$ in both norms. \second{The order drop for $N_v=12$ is due to the fact that the errors are evaluated with respect to a reference solution (not an exact one) computed on a grid with $N_v =13$.} Figure \ref{fig:test_staz_psi}, on the left, shows the approximation of the solution to the adjoint equation for $t=0,0.25,0.5,0.75$ and the terminal condition $t=1$, while on \second{the} right the second-order convergence to the reference solution is shown. 

\begin{comment}
	\begin{table}[ht]
		\centering\begin{minipage}{.8\textwidth}
			\caption{Errors and convergence orders for the Fokker--Planck equation at the final time $t=T$. In all simulations $\Delta t = \Delta v/2$.}
			\label{tab1}
			\begin{tabular*}{\textwidth}{@{\extracolsep{\fill}}*{6}{l}@{\extracolsep{\fill}}}
				\toprule%
				$\Delta v$&$\Delta t$&$E_2(f(\cdot,T))$&$p_2$&$E_\infty(f(\cdot,T))$&$p_\infty$\\
				\midrule
				${0.2}&${0.1}&${8.8025\times 10^{-1}&$/$&${9.2021\times 10^{-1}&$/$\\
				\midrule
				${0.1}&${0.05}&${1.8273\times 10^{-1}&$2.27$&${1.8447\times 10^{-1}&$2.32$\\
				\midrule
				${0.05}&${0.025}&${3.7548\times 10^{-2}&$2.28$&${3.3027\times 10^{-2}&$2.48$\\
				\midrule
				${0.025}&${0.0125}&${9.9255\times 10^{-3}&$1.92$&${9.6183\times 10^{-3}&$1.78$\\
				\bottomrule
		\end{tabular*}\end{minipage}
	\end{table}
\end{comment}
\begin{table}[ht]
	\centering\begin{minipage}{.8\textwidth}
		\caption{Errors and convergence orders for the Hamilton--Jacobi equation at the initial time $t=0$. In all simulations $\Delta t = \Delta v$.}
		\label{tab2}
		\begin{tabular*}{\textwidth}{@{\extracolsep{\fill}}*{6}{l}@{\extracolsep{\fill}}}
			\toprule%
			$N_v$&$E_2(\psi(\cdot,0))$&$p_2$&$E_\infty(\psi(\cdot,0))$&$p_\infty$ & \second{CPU time $(s)$}\\
			\midrule
			$4$&$8.8348\times 10^{-2}$&$/$&$8.8453\times 10^{-2}$&$/$ & $1.93 \times 10^{-3}$\\
			\midrule
			$5$&$2.1201\times 10^{-2}$&$2.05$&$2.0812\times 10^{-2}$&$2.09$ & $3.78 \times 10^{-3}$\\
			\midrule
			$6$&$4.8843\times 10^{-3}$&$2.12$&$5.0052\times 10^{-3}$&$2.06$& $1.22 \times 10^{-2}$\\
			\midrule
			$7$&$1.1746\times 10^{-3}$&$2.06$&$1.2131\times 10^{-3}$&$2.04$ & $4.74\times 10^{-2}$\\
			\midrule
			$8$&$2.8294\times 10^{-4}$&$2.05$&$2.9106\times 10^{-4}$&$2.06$ & $2.78\times 10^{-1}$\\
			\midrule
			$9$&$6.6038\times 10^{-5}$&$2.10$&$6.7338\times 10^{-5}$&$2.11$ & $3.19$\\
			\midrule
			$10$&$1.4541\times 10^{-5}$&$2.18$&$1.4293\times 10^{-5}$&$2.24$ & $3.74\times 10^{1}$\\
			\midrule
			$11$&$3.1253\times 10^{-6}$&$2.22$&$2.4059\times 10^{-6}$&$2.57$ & $4.69\times 10^{2}$\\
			\midrule
			$12$&$9.2072\times 10^{-7}$&$1.76$&$6.8016\times 10^{-7}$&$1.82$ & $5.85\times 10^3$\\
			\bottomrule
	\end{tabular*}\end{minipage}
\end{table}
\subsection{Numerical experiments in social dynamics}\label{sec:qualit}
\begin{comment}
	\begin{equation}\label{eq:opifp_old}
		\begin{cases}
			\partial_t f = \partial_v \mathcal{F}_v[f](v,c,t)+\partial_c \mathcal{F}_c[f](v,c,t) &(v,c,t) \in \Omega\times\mathbb{R}^2_+,\\
			f(v,c,0) = f_0(v,c), &(v,c) \in \Omega\times\mathbb{R}_+\\
			\left(\mathcal{V}[f,u]f + \nabla\left(D(v)f\right)\right)\cdot \vec{n} = 0, &x = \partial \Omega, c\in \mathbb{R}_+
		\end{cases}
	\end{equation}
	where the fluxes are defined as
	\begin{equation}\label{eq:fluxes_old}
		\begin{aligned}
			\mathcal{F}_v[f](v,c,t) &= -\left(\mathcal{E}[f](v,c,t) + u(v,c,t)\right) f + \frac{\sigma^2}{2}\partial_v\left(D(v)f\right),\\
			\mathcal{F}_c[f](v,c,t) &= \Phi(c)f + \frac{\nu^2}{2}\partial_c\left(c^2f\right),
		\end{aligned}
	\end{equation}
	
	and represents the density of agents who have an opinion $v\in[-1,1]=\Omega$ on a certain topic and the number of contacts on social media $c>0$ at time $t\geq 0$. As in \cite{albi24}, we suppose that the marginal distribution of $f$ with respect to contacts, at the stationary state is equal to a lognormal, meaning that we choose
	for a fixed choice of the parameters $\mu > 0$ and $\overline c > 0$.
\end{comment}
This concluding section is \second{focused on} two qualitative simulations in two spatial dimensions. In Section \ref{sec:opic}, the first spatial variable \second{is interpreted} as the opinion variable, while the second one represents the number of social media contacts of the agents in a population. The simulation shows how an external action on the opinions of certain groups of people can have a high impact in shaping the general sentiment of the entire population.

In Section \ref{sec:dueopi}, instead, the two variables $v_1$ and $v_2$ represent two opinions on different topics, one influencing the evolution of the other and vice versa. We show how an external force can modify the evolution of the general idea of the public and lead to consensus on an objective combination of opinions, even when the uncontrolled case offers a completely different outcome.

\second{\subsubsection{Comparison between the first-order and the second-order scheme}\label{sec:compo1o2}
	In this subsection, we compare the performance of a first-order and our second-order scheme, with the aim to show that the latter gives better results in terms of terminal cost compared to the first one.  The simulation is performed in one spatial dimension, with $\Omega=[-1,1]$. We define the quantities
	\[f_{min} = \min_{v,t}f(v,t), \quad \mbox{and}\quad E_{int} =\max_{t}\lvert \int_{\Omega} f(v,t)\mathrm{d} v - 1\rvert,
	\]
	which allow us to check the positivity of the structure-preserving scheme and the mass conservation.
	The problem we focused on is the following. In \eqref{eq:opifp} we set $P(v,v^*) = \chi(\lvert v-v^* \rvert \leq 0.1)$, $D(v) = 5\times 10^{-3}(1-v^2)$, and $f_0(v) = 0.5$, i.e. the initial condition is the uniform distribution in $\Omega$. In the cost \eqref{eq:cost_j}, $\gamma = 1$, $s(v) = 1$, $\overline v = 0.3$, and $T=4$. Tables \ref{tab3} and \ref{tab4} show the performance of the first-order and the second-order scheme, respectively. In these Tables, $N_v$ is such that the number of grid points of the space discretization is $2^{N_v}$, $\Delta x = 2/2^{N_v}$ and $\Delta t = \Delta x/10$. We can see that the second-order scheme reaches convergence in a smaller number of iterations, at the same value of the space step recovers a smaller optimal value of the cost functional $\mathcal{J}$, has a better conservation of the mass compared to the first-order one. Clearly, the main advantage of the first-order scheme is the speed: at the same space step size the execution time of each iteration is smaller.
}
\begin{table}[ht]
	\centering\begin{minipage}{.8\textwidth}
		\caption{Results of the first-order scheme. The table shows the computed optimal value of the cost functional $\mathcal{J}(\mathbf{u}^*;f_0)$, the number of iterations $k$, the CPU time per iteration, the minimum value of the function $f$ and the maximum error in the integral.}
		\label{tab3}
		\begin{tabular*}{\textwidth}{@{\extracolsep{\fill}}*{6}{l}@{\extracolsep{\fill}}}
			\toprule%
			$N_v$&$\mathcal{J}(\mathbf{u}^*;f_0)$&$k$&CPU time $(s)$&$f_{min}$ & $E_{int}$\\
			\midrule
			$4$&$2.95\times 10^{-1}$&$12$&$4.27\times 10^{-2}$&$1.48\times 10^{-12}$ & $2.44\times 10^{-15}$\\
			\midrule
			$5$&$2.61\times 10^{-1}$&$11$&$8.61\times 10^{-2}$&$4.04\times 10^{-25}$& $2.66\times 10^{-15}$\\
			\midrule
			$6$&$2.48\times 10^{-1}$&$9$&$2.78\times 10^{-1}$&$3.65\times 10^{-49}$& $5.11\times 10^{-15}$\\
			\midrule
			$7$&$2.43\times 10^{-1}$&$10$&$1.14$&$5.40\times 10^{-83}$& $2.02\times 10^{-14}$\\
			\bottomrule
	\end{tabular*}\end{minipage}
\end{table}

\begin{table}[ht]
	\centering\begin{minipage}{.8\textwidth}
		\caption{Results of the second-order scheme. The table shows the computed optimal value of the cost functional $\mathcal{J}(\mathbf{u}^*;f_0)$, the number of iterations $k$, the CPU time per iteration, the minimum value of the function $f$ and the maximum error in the integral.}
		\label{tab4}
		\begin{tabular*}{\textwidth}{@{\extracolsep{\fill}}*{6}{l}@{\extracolsep{\fill}}}
			\toprule%
			$N_v$&$\mathcal{J}(\mathbf{u}^*;f_0)$&$k$&CPU time $(s)$&$f_{min}$ & $E_{int}$\\
			\midrule
			$4$&$2.92\times 10^{-1}$&$7$&$1.09\times 10^{-1}$&$9.41\times 10^{-13}$ & $6.66\times 10^{-16}$\\
			\midrule
			$5$&$2.59\times 10^{-1}$&$8$&$2.14\times 10^{-1}$&$4.15\times 10^{-27}$& $8.88\times 10^{-16}$\\
			\midrule
			$6$&$2.47\times 10^{-1}$&$8$&$6.73\times 10^{-1}$&$3.52\times 10^{-54}$& $1.11\times 10^{-15}$\\
			\midrule
			$7$&$2.42\times 10^{-1}$&$8$&$2.69$&$0$& $8.88\times 10^{-16}$\\
			\bottomrule
	\end{tabular*}\end{minipage}
\end{table}

\subsubsection{Evolving opinions and social media contacts}\label{sec:opic}
\begin{figure}[h!]
	\centering
	\includegraphics[width=0.31\linewidth]{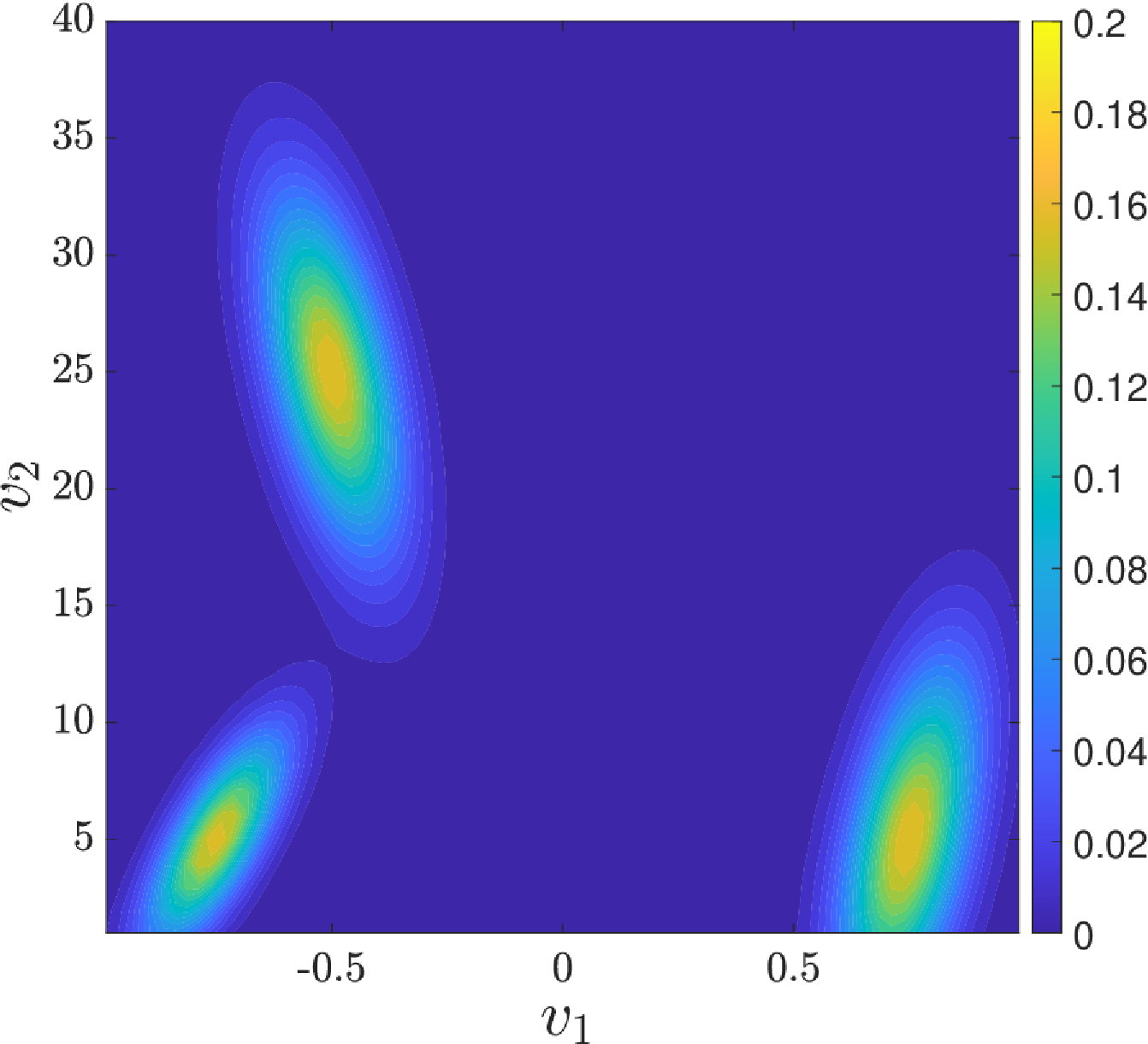}\quad \includegraphics[width=0.31\linewidth]{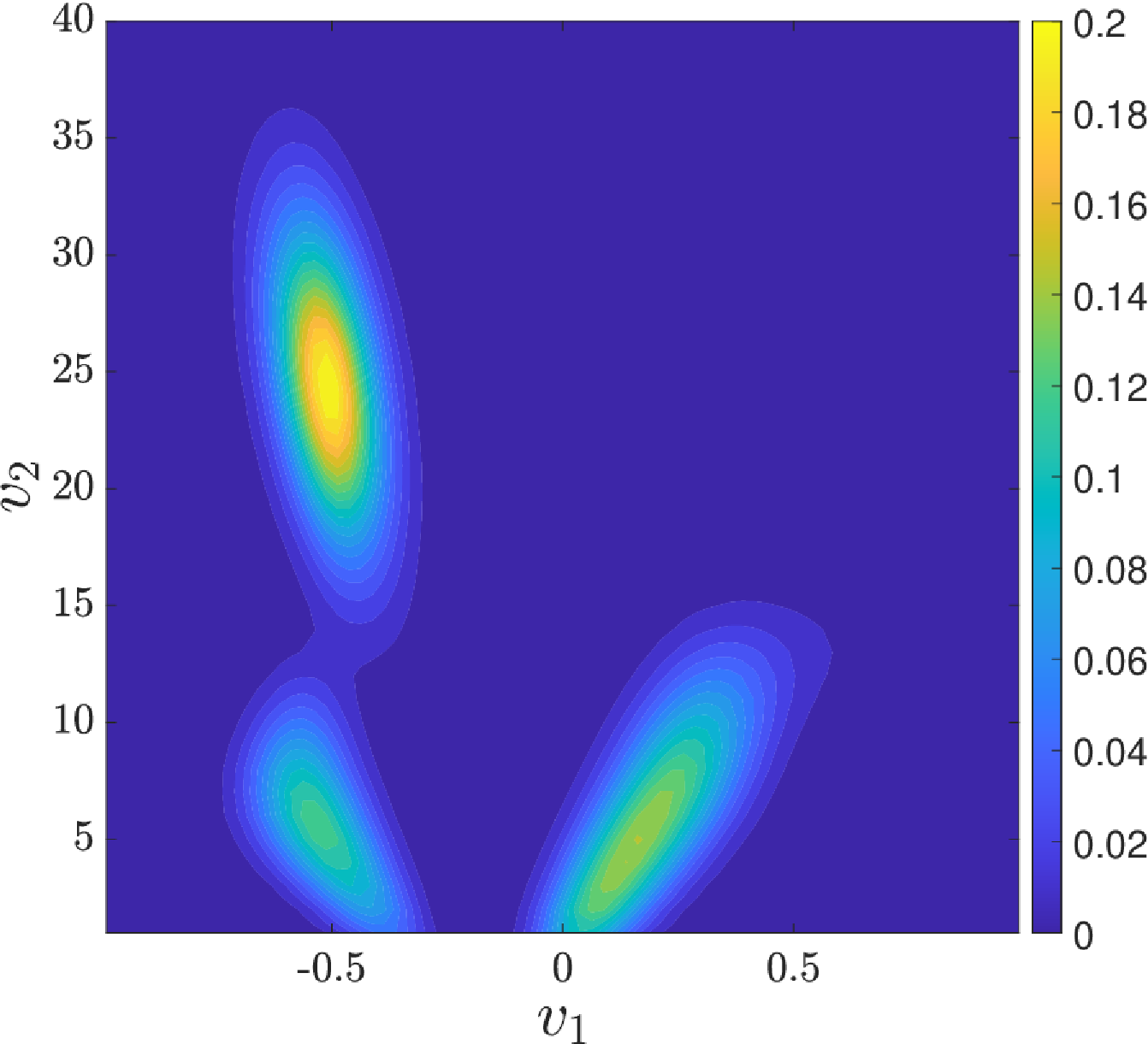}\quad \includegraphics[width=0.31\linewidth]{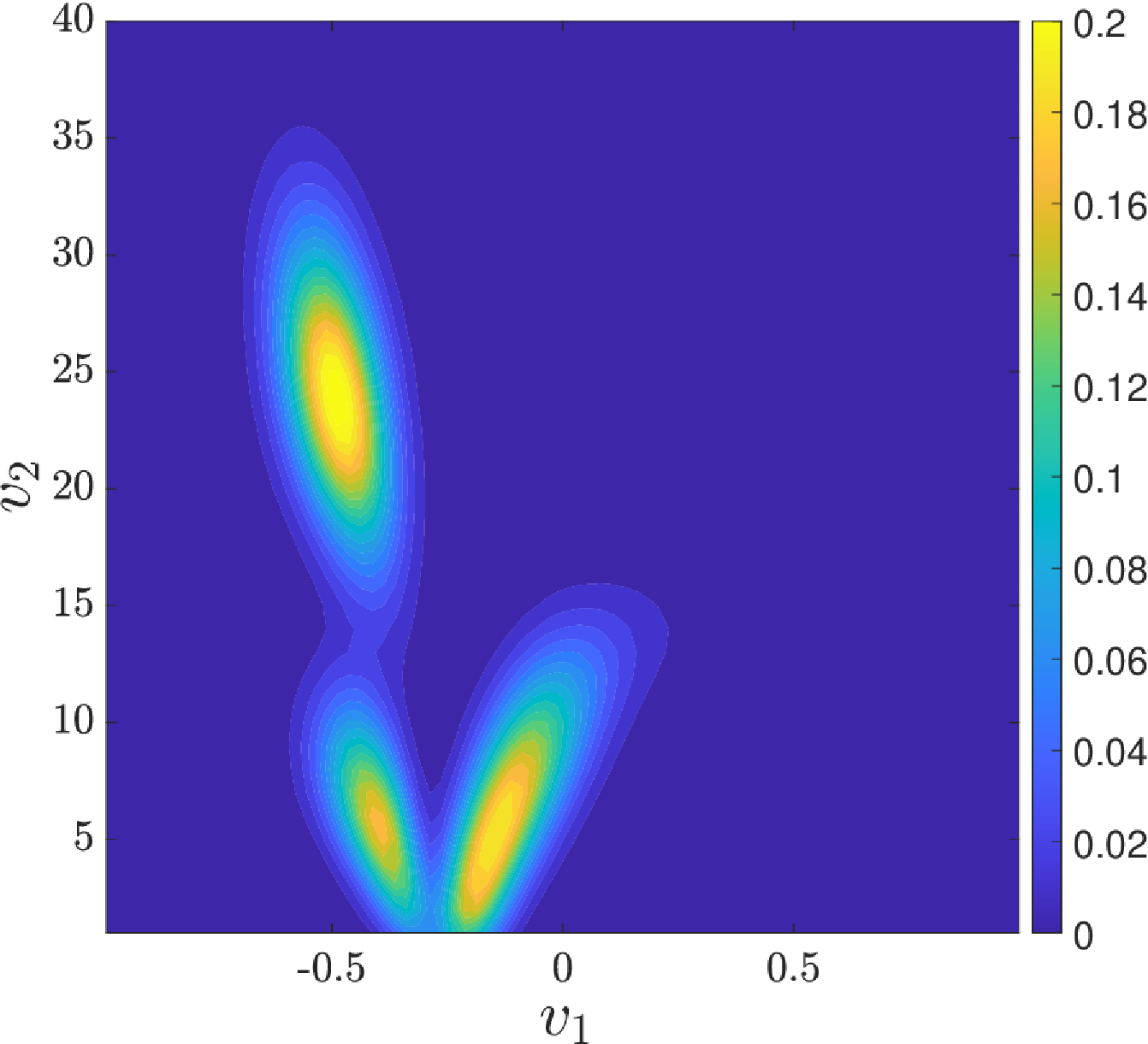}
	\caption{Uncontrolled density $f(\mathbf{v},t)$ at time $\second{t}=0$ (left), $t=1.5$ (center), and $t=3$ (right).}
	\label{fig:init_cond}
\end{figure}

\begin{figure}[h!]
	\centering
	\includegraphics[width=0.31\textwidth]{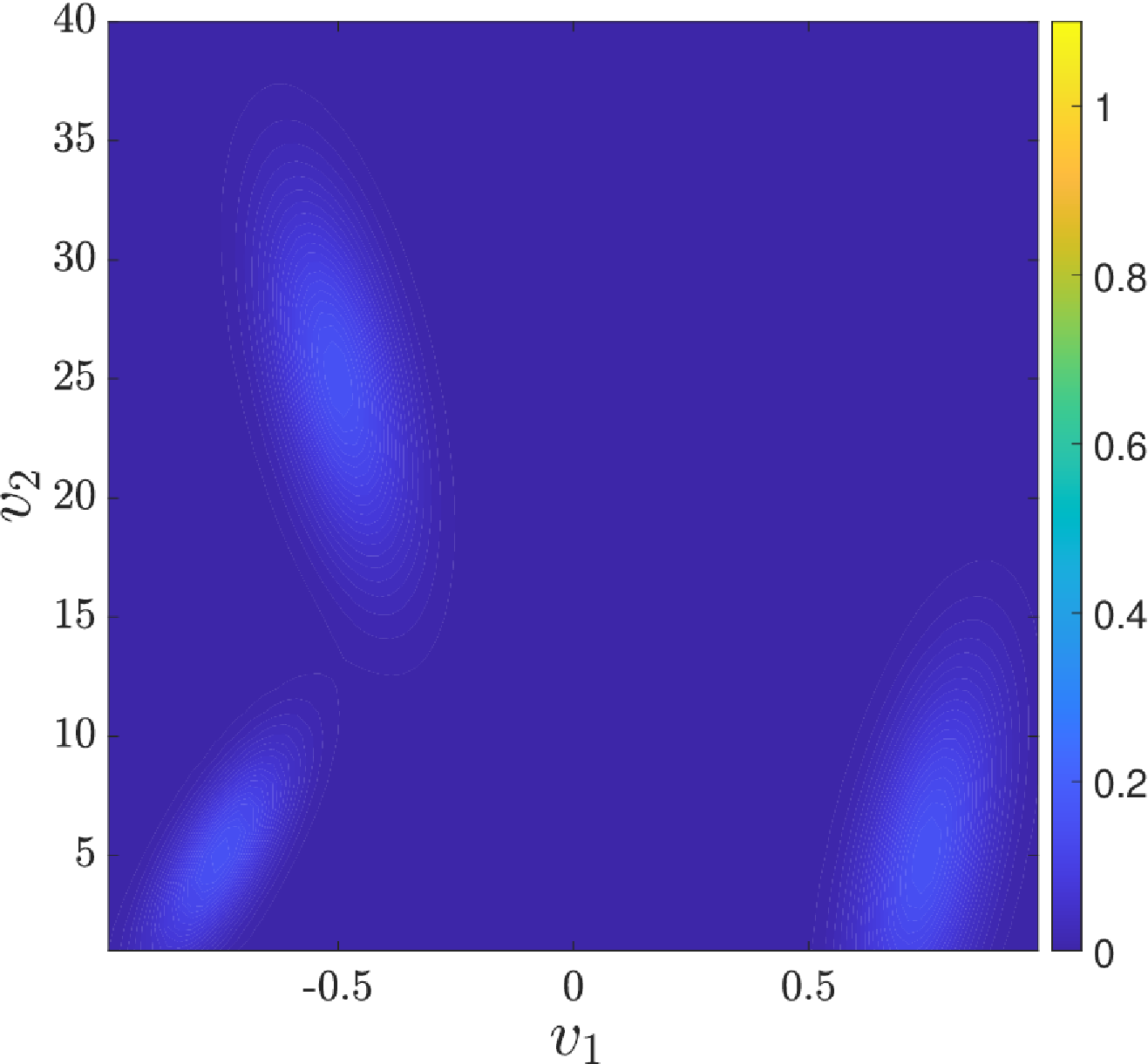}\quad
	\includegraphics[width=0.31\textwidth]{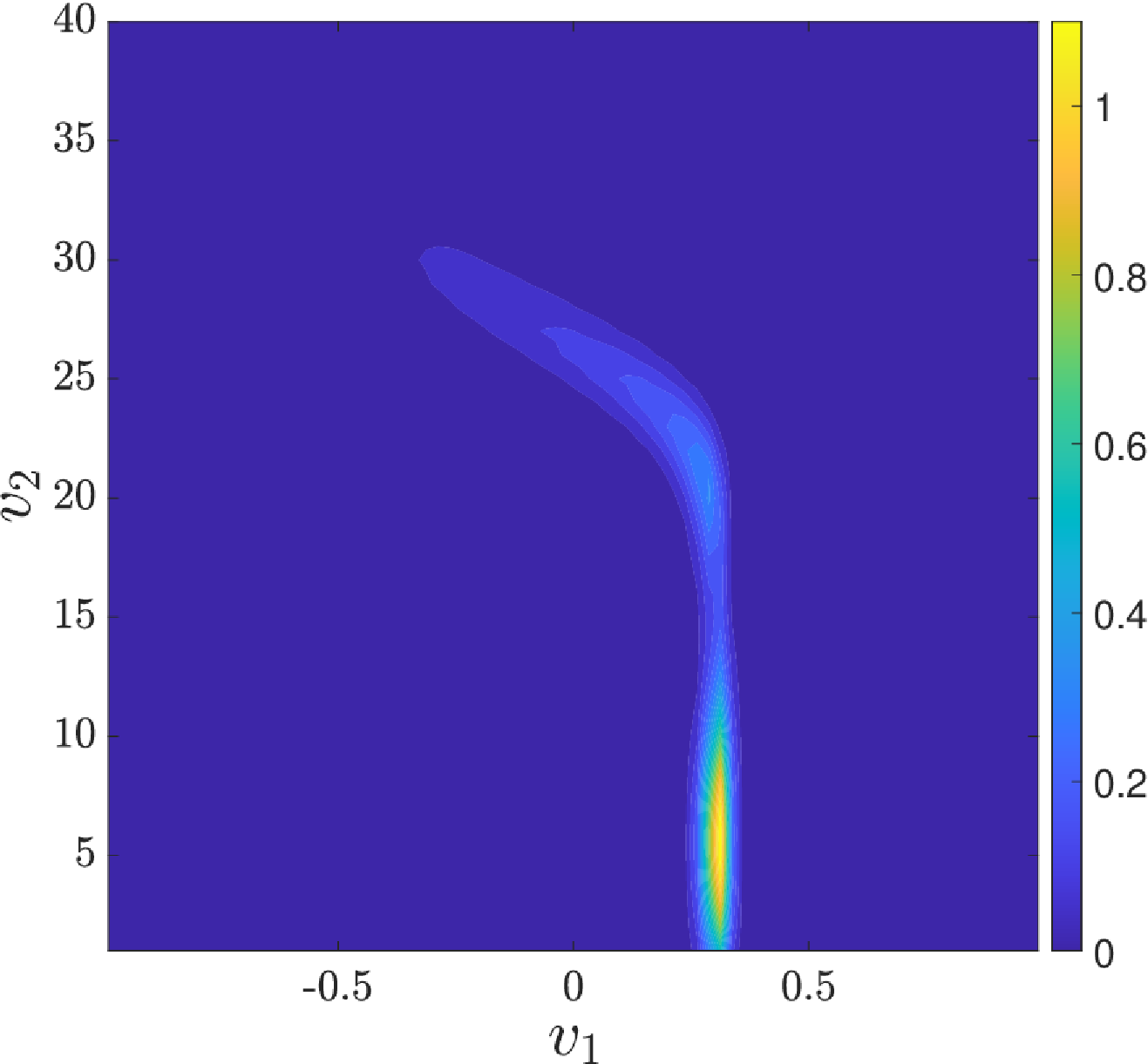}\quad
	\includegraphics[width=0.31\textwidth]{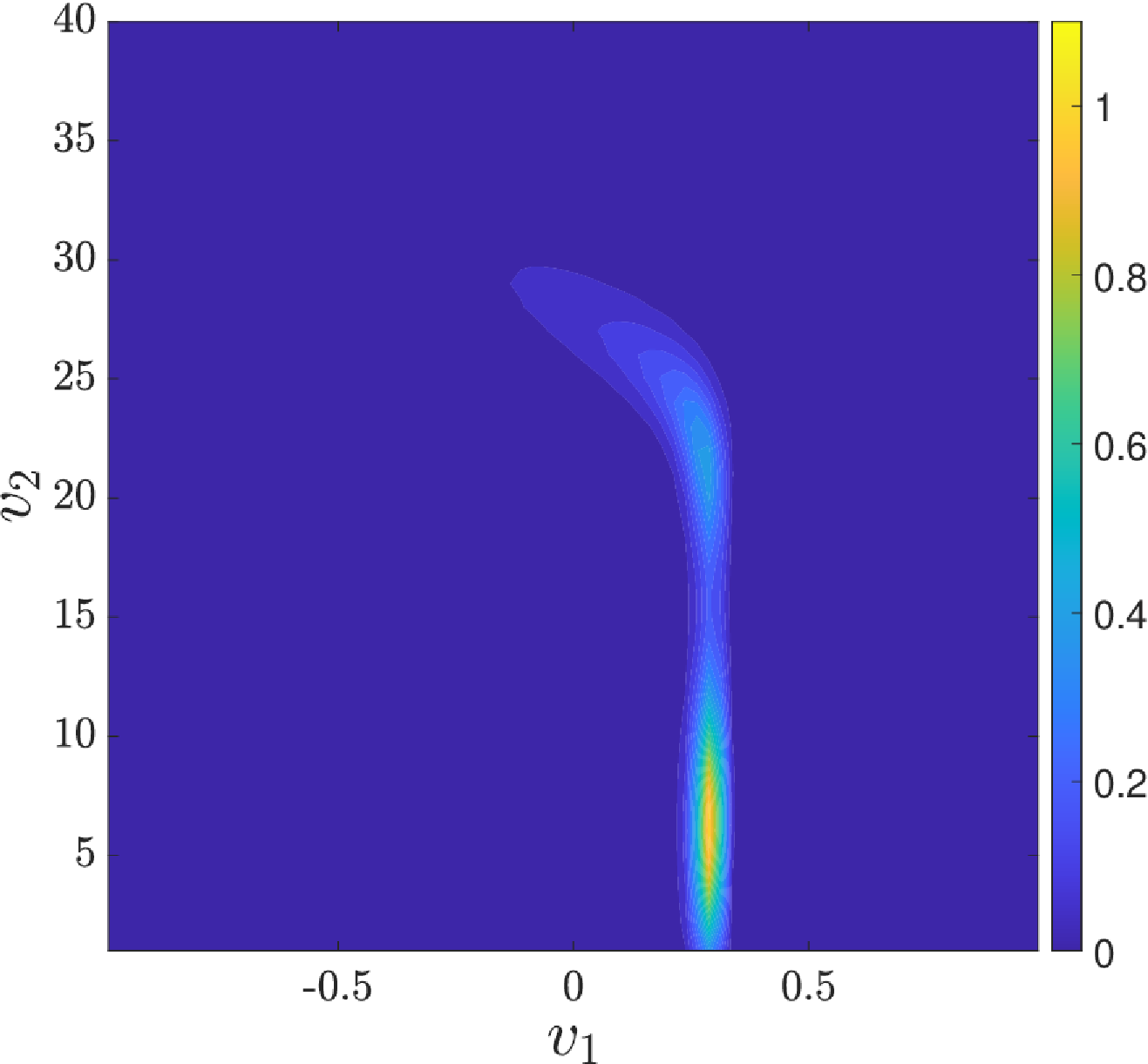}
	\\
	\centering
	\includegraphics[width=0.31\textwidth]{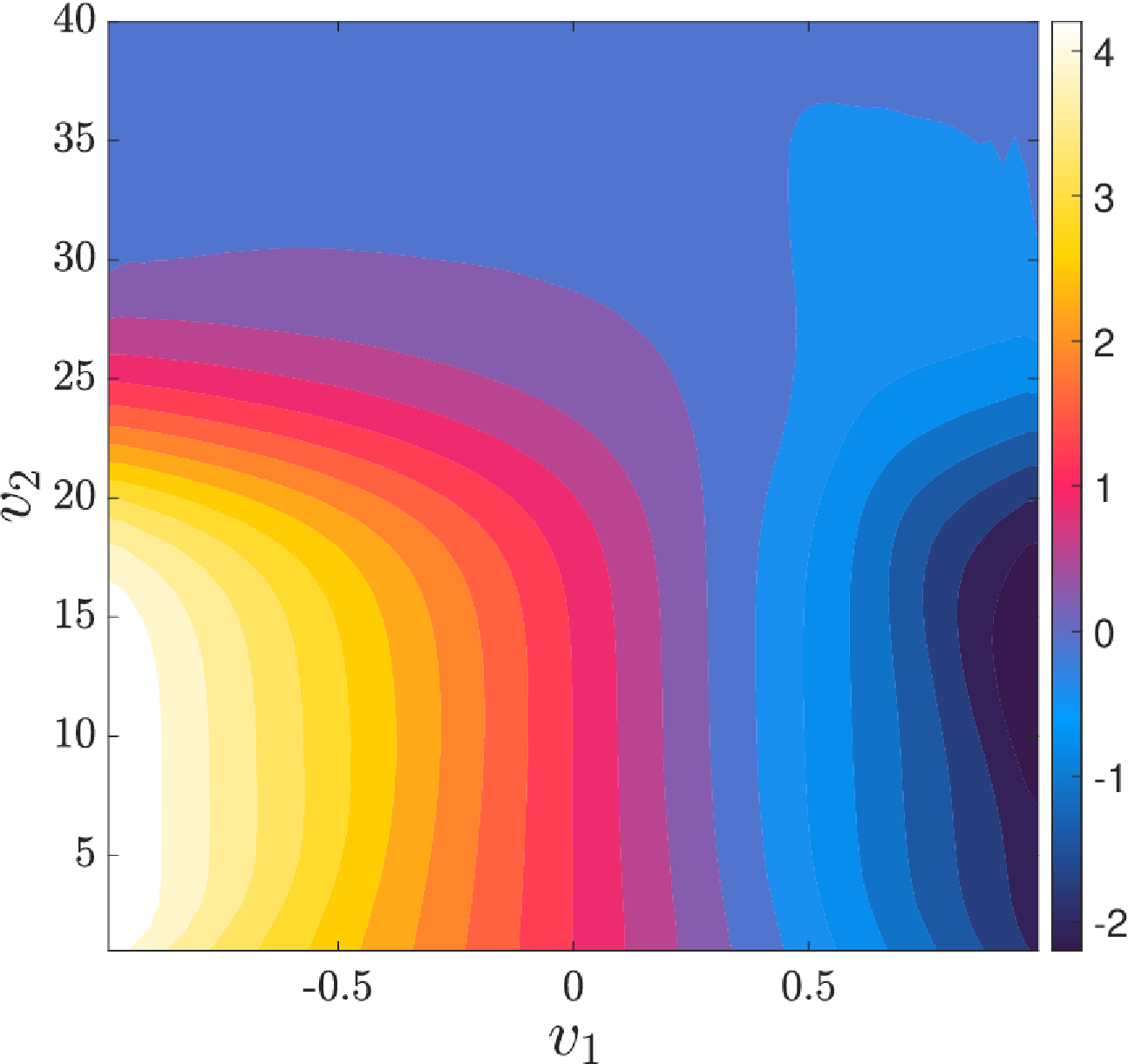}\quad
	\includegraphics[width=0.31\textwidth]{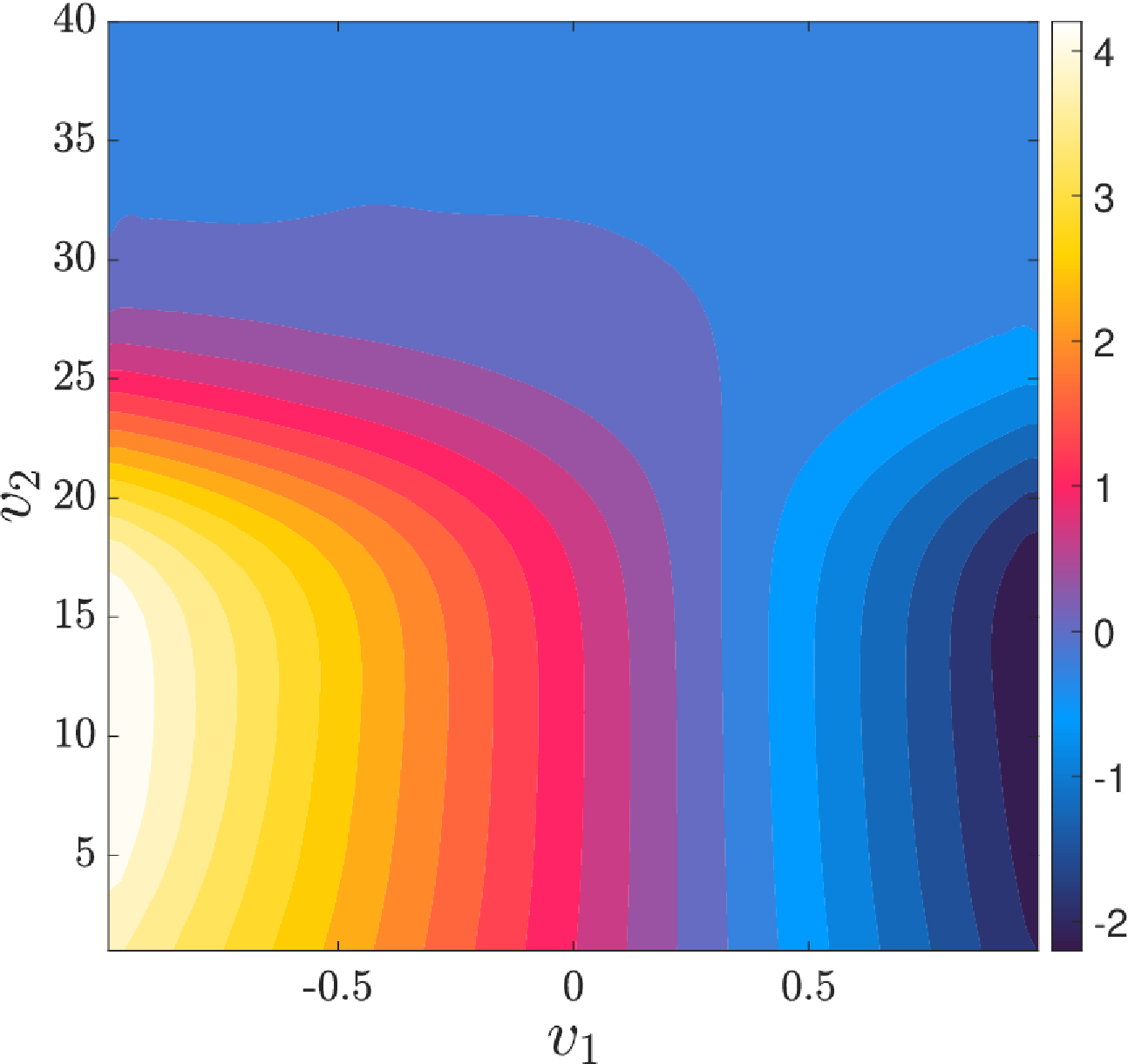}\quad
	\includegraphics[width=0.31\textwidth]{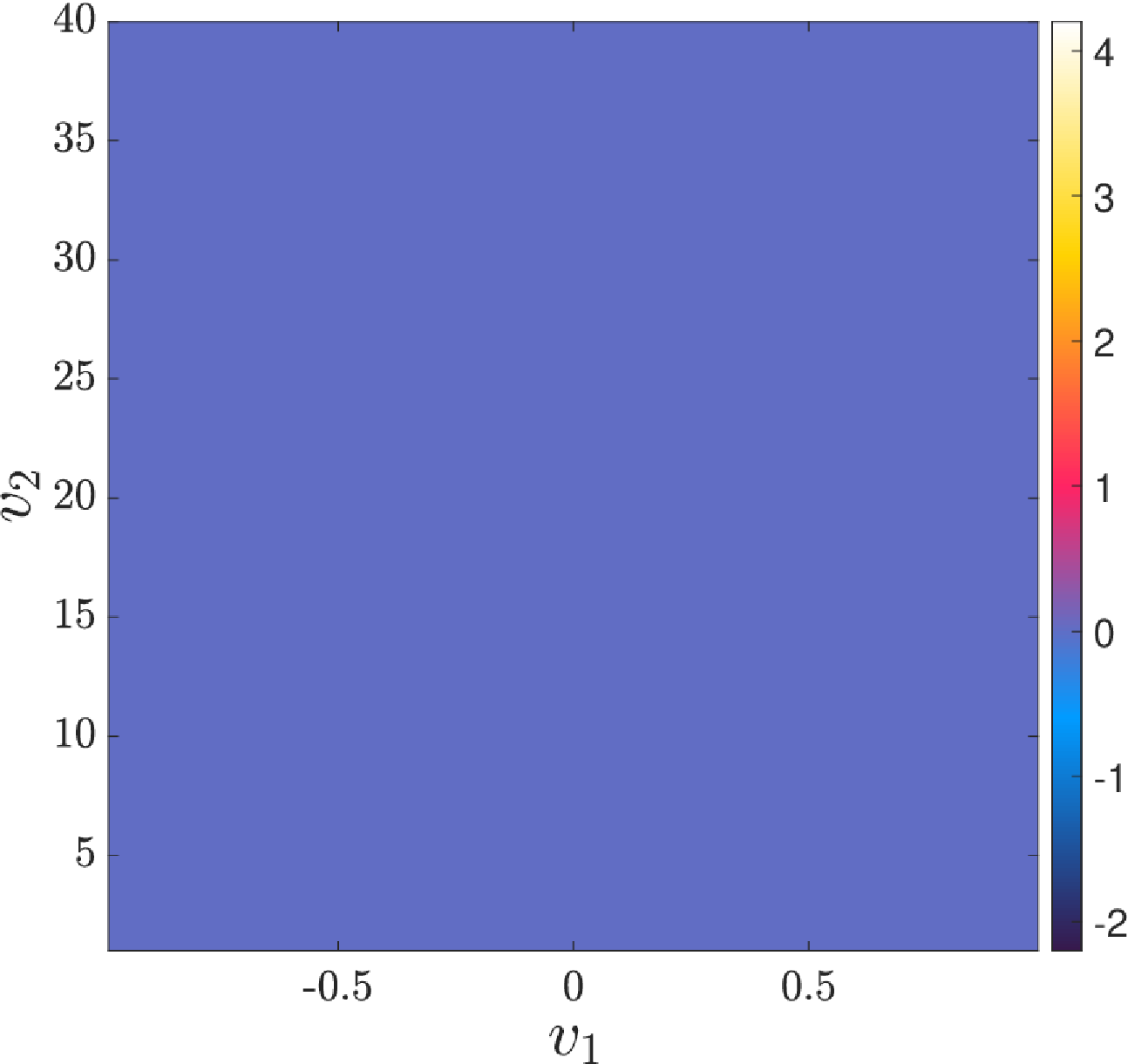}
	\caption{Case with $s_1(v_1,v_2)$ as in \eqref{eq:noinflcontr}. Contour plots of the density $f(\mathbf{v},t)$ (upper row), and the control $\mathbf{u}(\mathbf{v},t)$ (bottom row) for time $t=0$ (left), $t=1.5$ (center) and $t=3$ (right).}
	\label{fig:test3_1b}
\end{figure}

\begin{figure}[h!]
	\centering\includegraphics[width=0.31\textwidth]{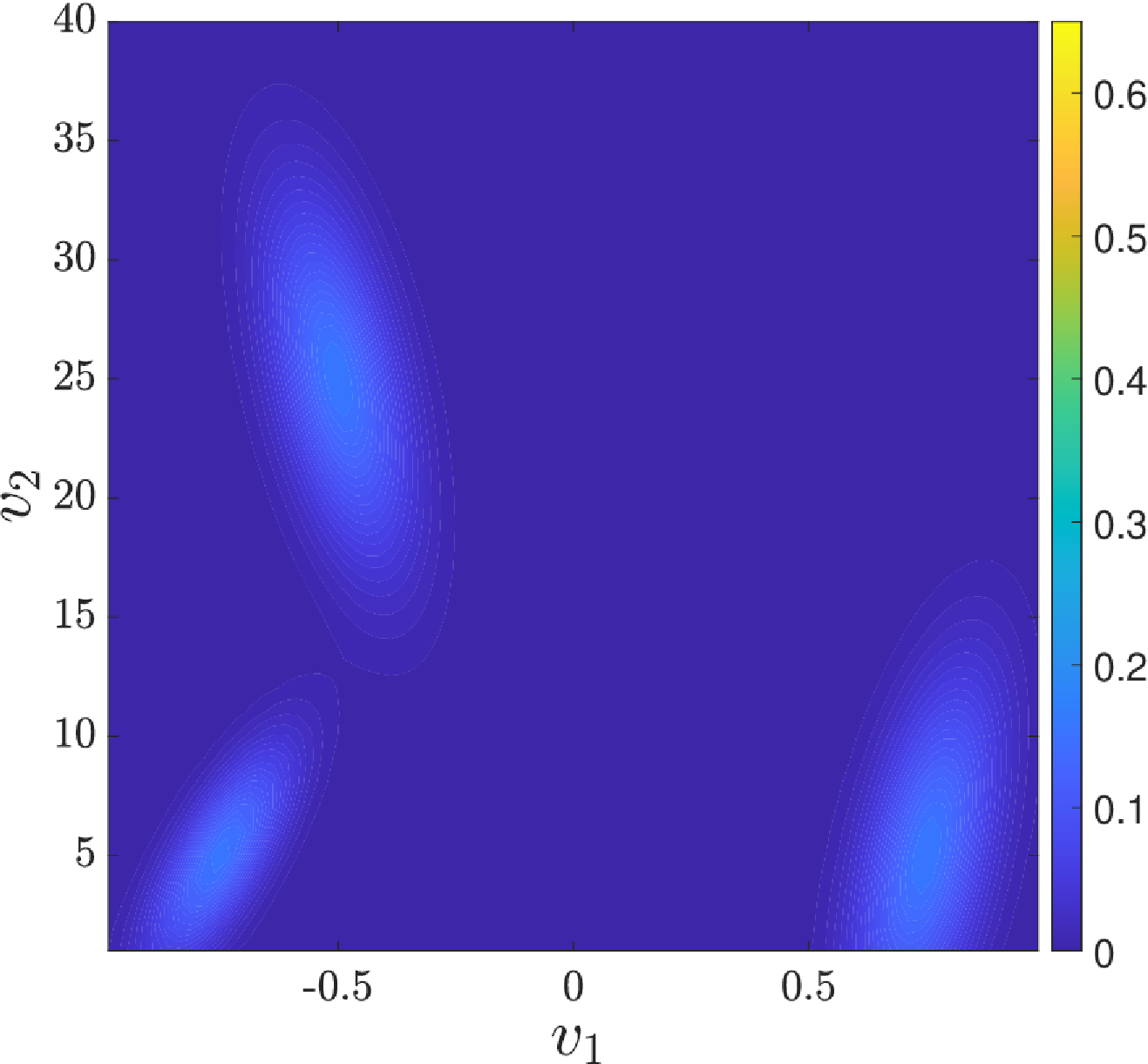}
	\quad	\includegraphics[width=0.31\textwidth]{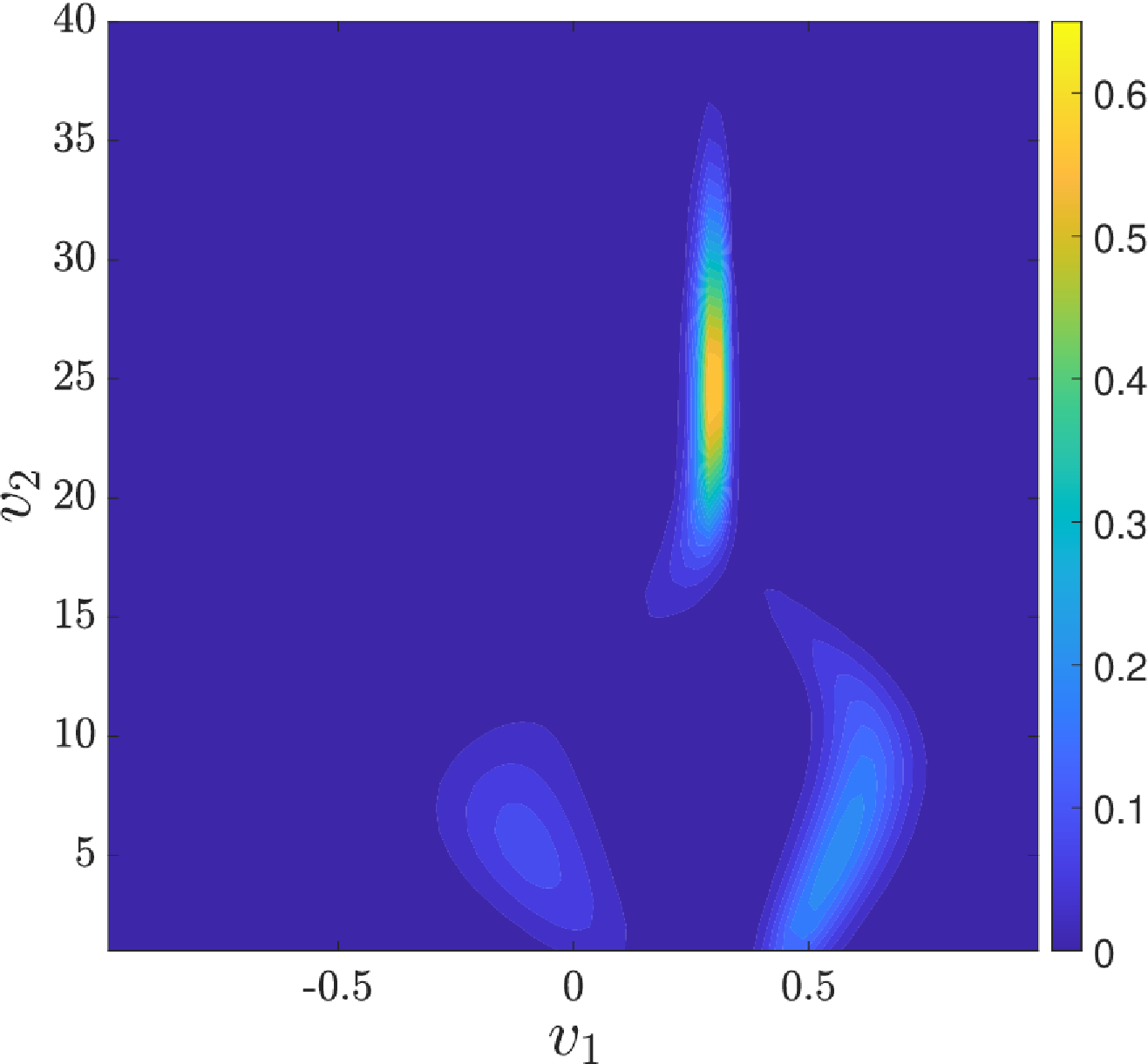}\quad \includegraphics[width=0.31\textwidth]{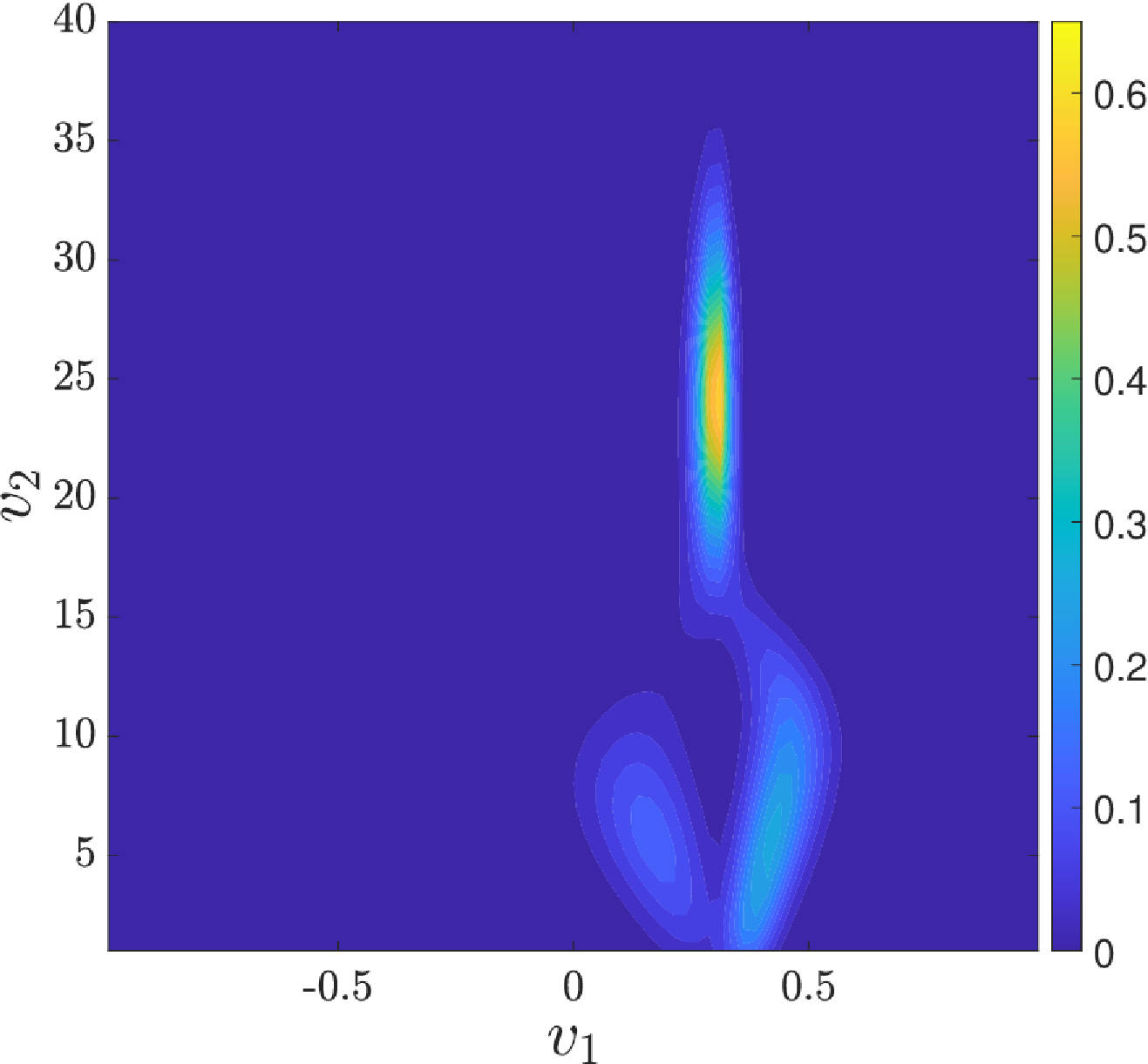}
	\\
	\centering
	\includegraphics[width=0.31\textwidth]{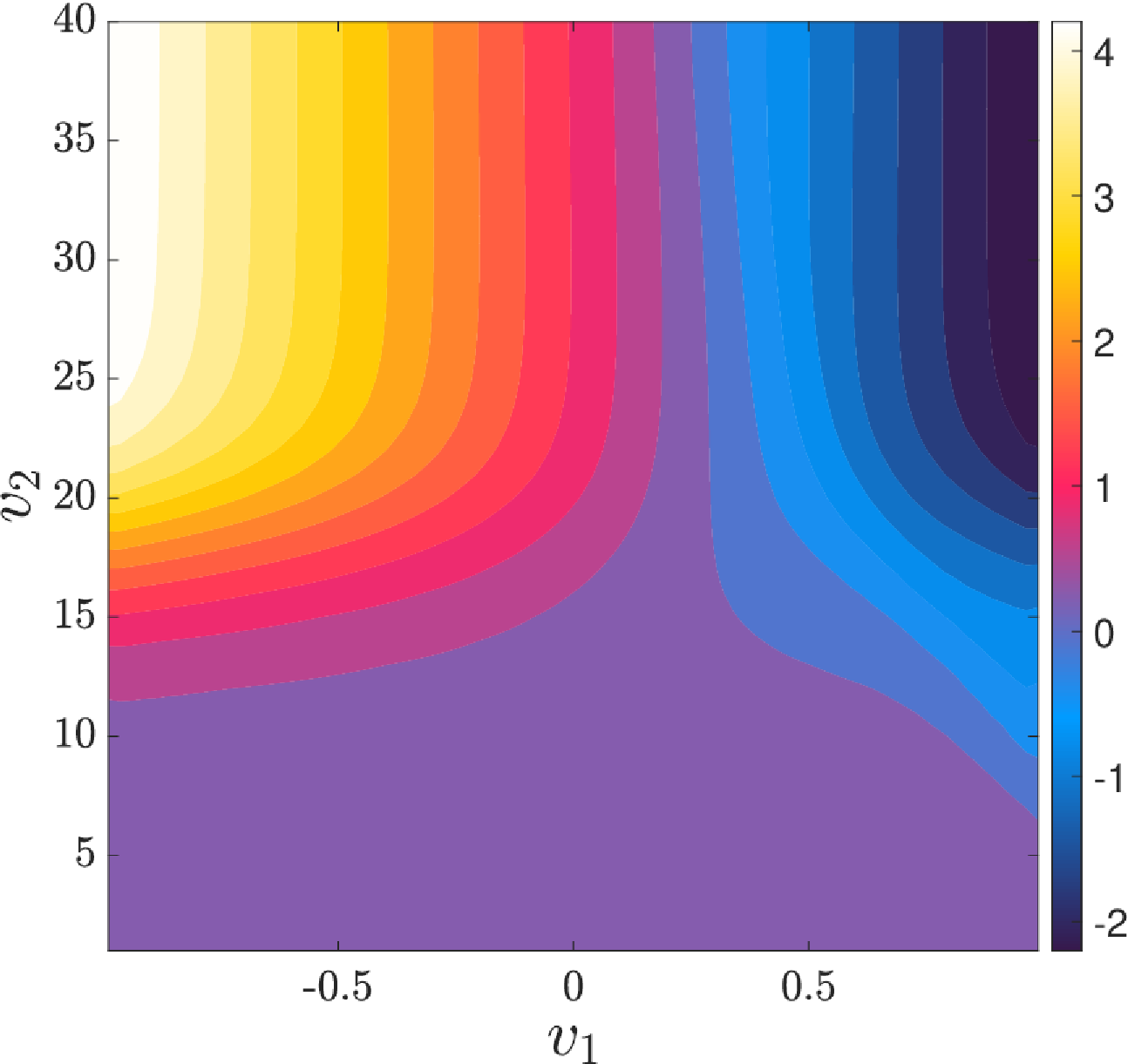}\quad\includegraphics[width=0.31\textwidth]{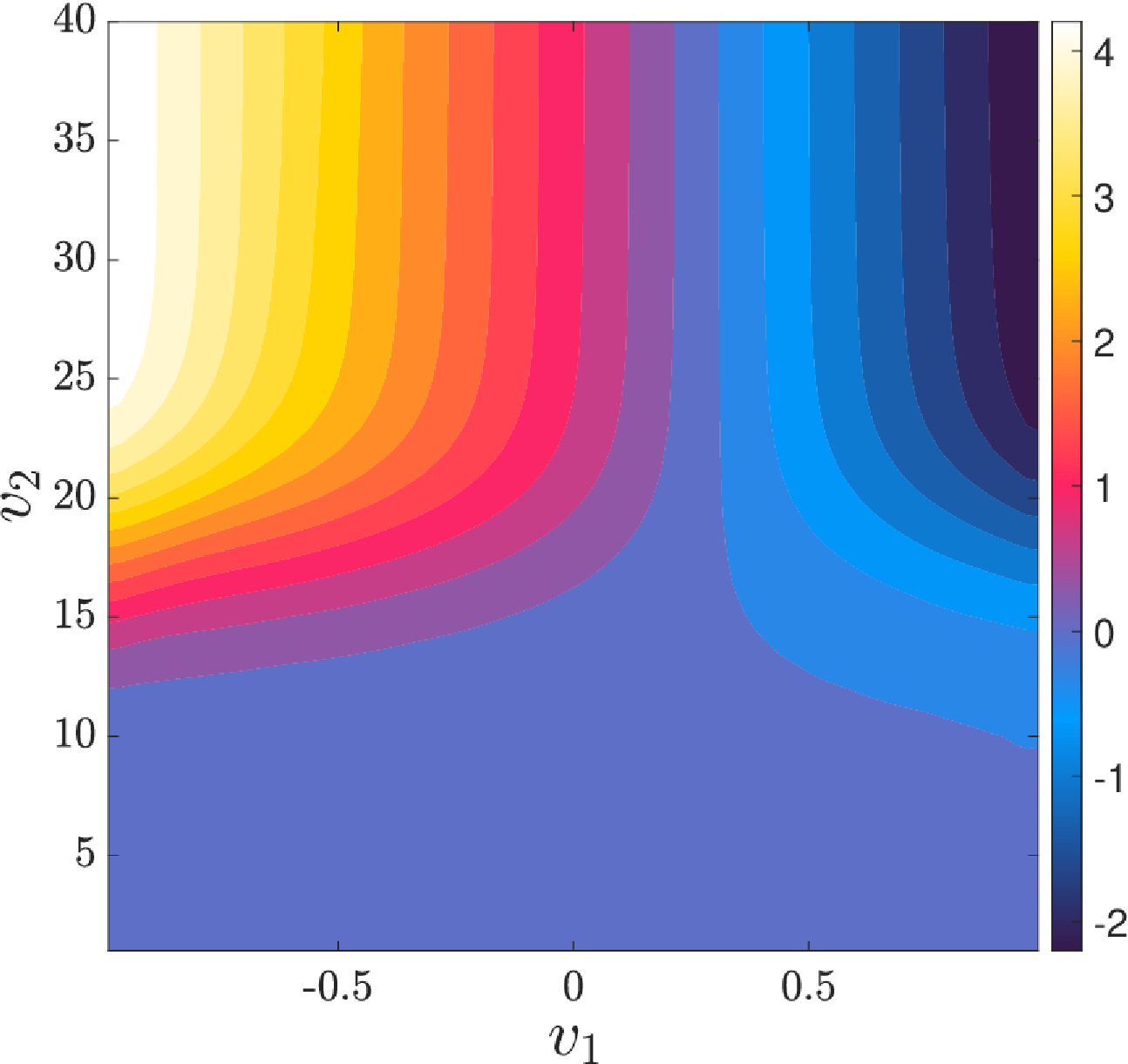}\quad\includegraphics[width=0.31\textwidth]{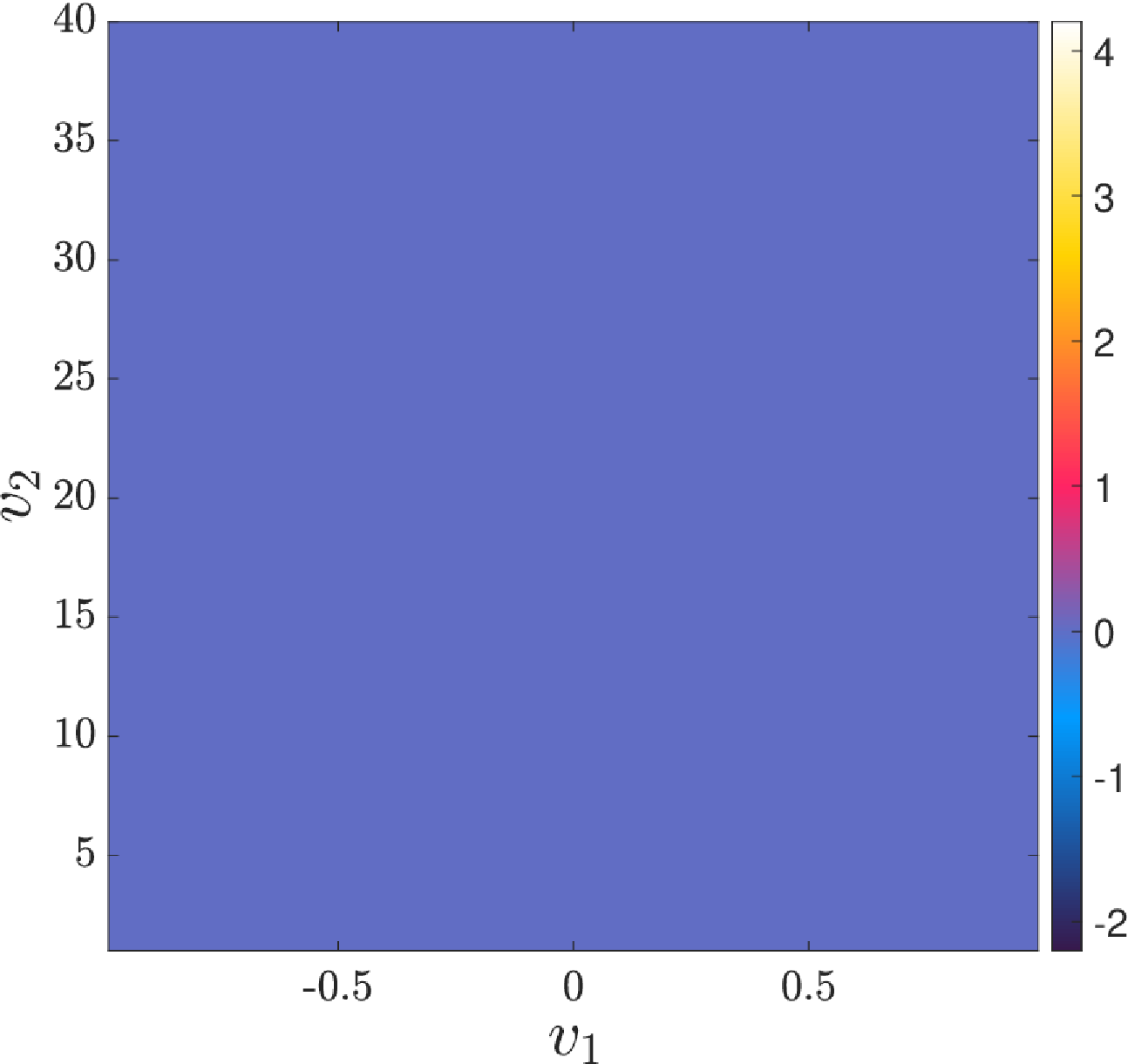}
	\caption{Case with $s_1(v_1,v_2)$ as in \eqref{eq:inflcontr}.  Contour plots of the density $f(\mathbf{v},t)$ (upper row), and the control $\mathbf{u}(\mathbf{v},t)$ (bottom row) at time $t=0$ (left), $t=1.5$ (center) and $t=3$ (right).}
	\label{fig:test3_2b}
\end{figure}

\begin{figure}[h!]
	\centering
	\includegraphics[width=0.45\textwidth]{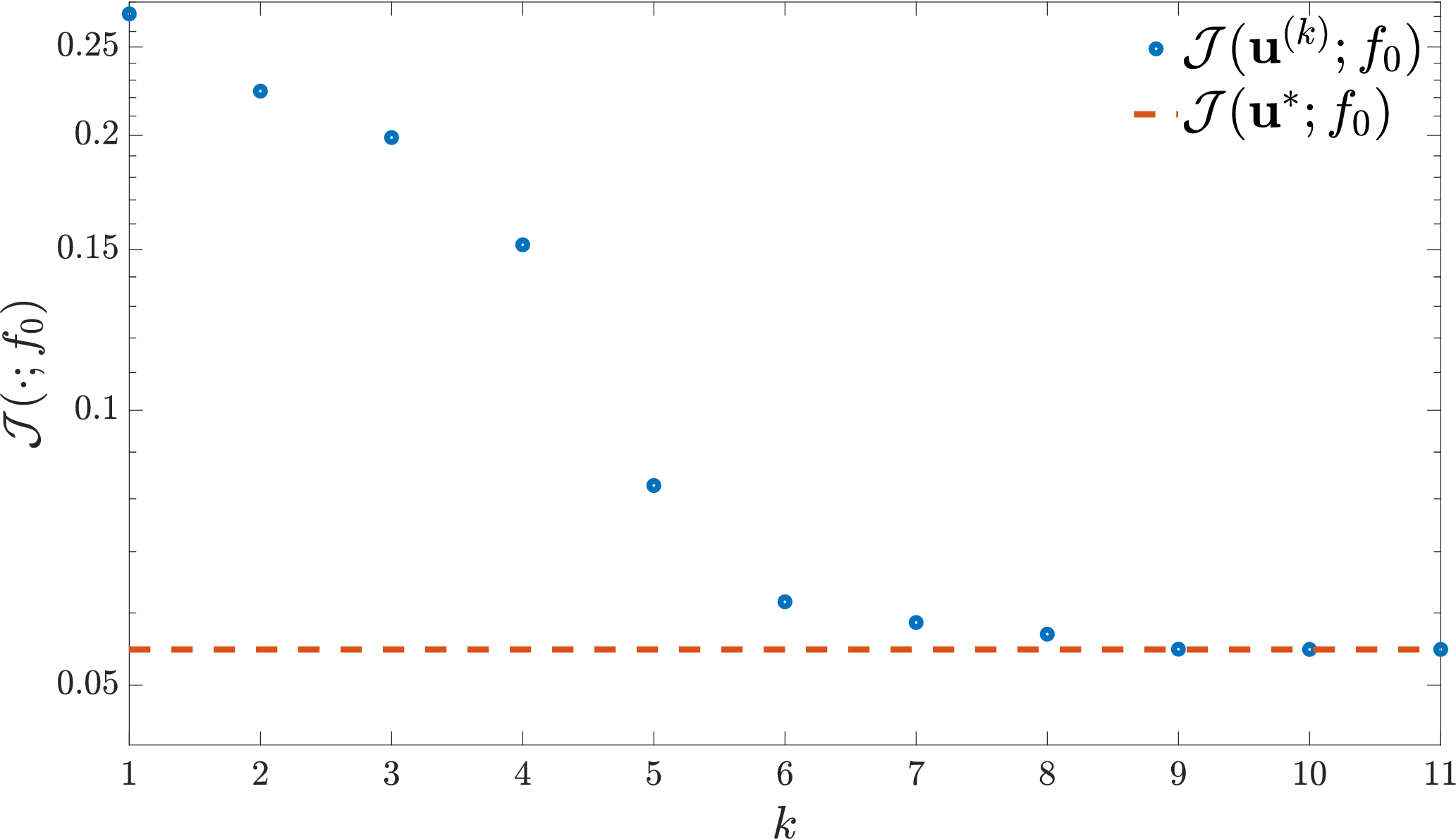}\quad
	\includegraphics[width=0.45\textwidth]{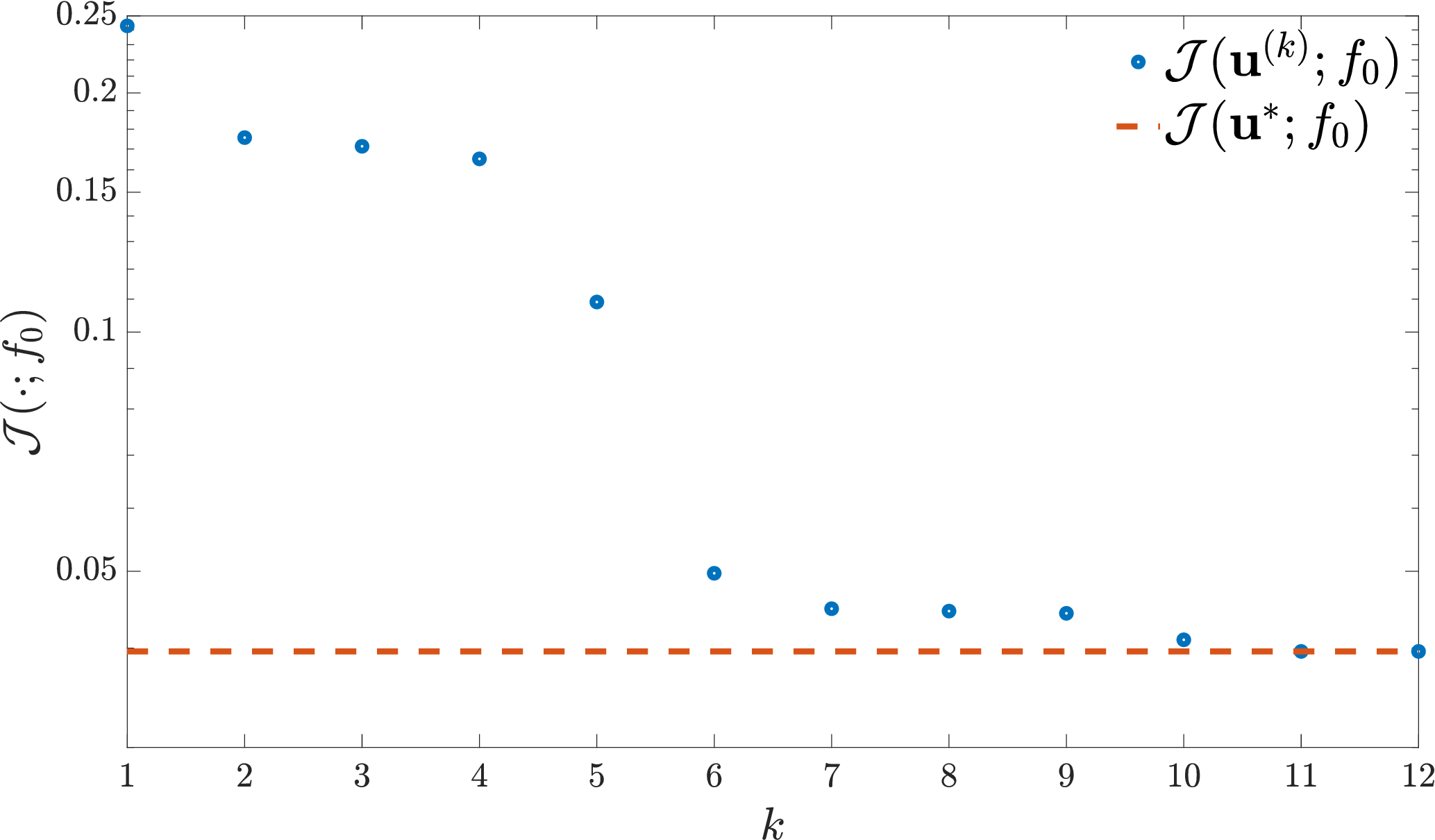}
	\caption{\second{Values of the cost functional $\mathcal{J}(\mathbf{u}^{(k)};f_0)$ at each iteration $k$ compared to the computed optimal value $\mathcal{J}(\mathbf{u}^*;f_0)$ for \eqref{eq:noinflcontr} (left), and \eqref{eq:inflcontr} (right). On the left we have $\mathcal{J}(\mathbf{u}^*;f_0) = 5.47 \times 10^{-2}$, while on the right  $\mathcal{J}(\mathbf{u}^*;f_0) = 3.96\times 10^{-2}$.}}
	\label{fig:test_gaussiane_2}
\end{figure}

In this simulation\second{,} we suppose that $\Omega_1=[-1,1]$ \second{and} $\Omega_2 = [1,v_{2,M}]$ with $v_{2,M} = 40$. Heuristically, \second{following the construction in \cite{albi24}}, we can consider the first variable $v_1$ to be the population opinion on a matter and the second one $v_2$ to be the number of social media contacts.
%we define the drift vector as
We \second{set} in \eqref{eq:drift_general} $h_1(v_1,v_2)=0$ and 
\[
P_{\second{1}}(v_1,v_2,v_1^*,v_2^*) = \chi\left(\lvert v_1-v_1^*\rvert \leq \Delta \frac{v_2^*}{v_2+v_2^*}\right)\frac{v_2^*}{v_2+v_2^*},
\]
\second{meaning that the interactions among the population happen if the distance between the opinions of the agents is below a certain threshold $\Delta v_2^*/(v_2+v_2^*)$. Such threshold is bigger if the interaction happens with an agent with a higher number of followers, conversely it is smaller if the interacting agent has a lower number of contacts, modeling the fact that people on the internet are more prone to adjust their opinions towards the ideas claimed by influent agents, while they are less eager to compromise with users with a smaller network. Then we set}
$P_2(v_1,v_2,v_1^*,v_2^*)=0$, $u_2(v_1,v_2,t) = 0$, and 
\begin{equation}\label{eq:drift_3_2}
	h_2(v_1,v_2) = \frac{\mu}{2}\log\left(\frac{v_2}{\hat v_2} \right) v_2, 
\end{equation}
where we set $\mu = 0.1$ and $\hat v_2 = 20$, and diffusion matrix
\[
\mathcal{D}(v_1,v_2) =
\begin{bmatrix}
	\frac{\sigma^2}{2}\left(1-v_1^2\right) & 0 \\ 0 & \frac{\nu^2}{2}v_2^2
\end{bmatrix},
\]
with $\Delta = 2$, $\sigma^2 = \nu^2 = 2\times 10^{-3}$.
Notice that this choice of $\mathcal{G}[f;\mathbf{u}](v_1,v_2,t)$ leads to $\mathcal{Q}_2[f](v_1,v_2,t)=0$ in \eqref{eq:hj}. \second{The choice of $h_2$ made in \eqref{eq:drift_3_2} is such that, at the stationary state, the distribution of connections on the network is a log-normal, which is coherent with the data obtained from $\mathbb{X}$ (even though we selected different parameters with respect to the actual ones of the social network). For further details on the topic we refer the reader to \cite{albi24}.} The initial condition is the sum of three bivariate gaussians
\begin{gather*}
	f_0(v_1,v_2) =K_0\left(e^{-\frac{1}{2(1-0.5^2)}\left(\frac{(v_1+0.5).^2}{0.01} + \left(\frac{v_1+0.5}{0.1}\right)\left(\frac{v_2-10}{5}\right)
		+ \frac{(v_2-10)^2}{25}\right)} \right.\\
	+ e^{-\frac{1}{2(1-0.5^2)}\left(\frac{(v_1-0.75)^2}{0.01)} - \left(\frac{v_1-0.75}{0.1)}\right)\left(\frac{v_2-50}{5}\right) +\frac{(v_2-50)^2}{25}\right)} \\
	\left.+ e^{-\frac{1}{2(1-0.75^2)}\left(\frac{(v_1+0.75)^2}{0.01} -1.5\left(\frac{v_1+0.75}{0.1}\right)\left(\frac{v_2-50}{5}\right) + \frac{(v_2-50)^2}{25}\right)}\right),
\end{gather*}
The constant $K_0$ normalizes the initial total mass, i.e. 
\[\int_{-1}^1\int_1^{v_{2,M}}f_0(v_1,v_2)\mathrm{d}v_1\mathrm{d}v_2 = 1.\]
We fix the final time $T=3$ \first{and $\Delta t = \Delta v_1\Delta v_2/(10(\Delta v_1+\Delta v_2))$}. Figure \ref{fig:init_cond} shows the initial condition $f_0$ (left), the numerical approximation of the uncontrolled density at time  $t=1.5$ (center) and $t=3$ (right). With this choice of interaction function and in the absence of control, the density tends to concentrate mainly around two negative opinions. 
In the functional \eqref{eq:cost_j}\second{,} we \second{set} $\overline v_1 = 0.3$ and $\gamma = 0.05$. 
We fix $s_2(v_1,v_2)=0$ and we perform simulations for two different values of the function $s_1(v_1,v_2)$: in the first case 
\begin{equation}\label{eq:noinflcontr}
	s_1(v_1,v_2) = \frac{1}{1+e^{-\frac 12 (\hat v_2 - v_2)}},
\end{equation}
meaning that the control acts mainly on the agents with less than $\hat v_2$ contacts, while in the second case we set
\begin{equation}\label{eq:inflcontr}
	s_1(v_1,v_2) = \frac{1}{1+e^{-\frac 12 (v_2 - \hat v_2)}},
\end{equation}
so the control acts more intensely on agents with more than $\hat v_2$ contacts. Figure \ref{fig:test3_1b} shows a comparison between the controlled density (upper row) and the control function (bottom row) at time $t=0$ (left), $t=1.5$ (center) and $t=3$ (right) for the choice in \eqref{eq:noinflcontr}. We can see how the control tends to concentrate the agents with a lower number of contacts around $v_1=0.3$, as expected, but part of the influent population (with a higher value of $v_2$) stays on the negative half of $\Omega_1$. Figure \first{\ref{fig:test3_2b}} shows a comparison between the controlled density (upper row) and the control function (bottom row) at time $t=0$ (left), $t=1.5$ (center) and $t=3$ (right) for the choice in \eqref{eq:inflcontr}. Again, we see how the control is able to concentrate the density around $v_1=0.3$, acting primarily on the part of the population with a higher number of contacts. Due to the expression of the interaction function $P$, also the segment of the population with less contacts gets influenced and modifies its opinion, ending up with a positive value of $v_1$. \first{It is worth noticing that the terminal condition $\mathbf{u}(\mathbf{v},T) = 0$ follows from equation \eqref{eq:optimcond} and from the absence of a terminal cost in \eqref{eq:cost_j}}. Figure \ref{fig:test_gaussiane_2} shows the value of the cost functional $\mathcal{J}(\mathbf{u}^{(k)};f_0)$ for the various iterations $k$ of the gradient descent method compared to the computed optimal value $\mathcal{J}(\mathbf{u}^{*};f_0)$, both for the choice of $s_1$ in \eqref{eq:noinflcontr} (left) and in \eqref{eq:inflcontr} (right). \second{We can see how the computed optimal cost for \eqref{eq:noinflcontr} is higher than the one for \eqref{eq:inflcontr} and this is in accordance with the heuristic interpretation of the two choices of $s_1(v_1,v_2)$. In fact, it is more efficient (i.e. \emph{cheaper}) to control the opinion of the most popular agents if we want to steer the public opinion towards a certain value than try to obtain the same results controlling the less influent part of the population.}
\begin{comment}
	\begin{figure}[h!]
		\centering
		\includegraphics[width=0.31\textwidth]{testnoinfl_T1.eps}\quad
		\includegraphics[width=0.31\textwidth]{testnoinfl_T2.eps}\quad
		\includegraphics[width=0.31\textwidth]{testnoinfl_T3.eps}
		\\
		\centering
		\includegraphics[width=0.31\textwidth]{noinfl_contr_T1.eps}\quad
		\includegraphics[width=0.31\textwidth]{noinfl_contr_T2.eps}\quad
		\includegraphics[width=0.31\textwidth]{noinfl_contr_T3.eps}
		\caption{Case with $s_1(v_1,v_2)$ as in \eqref{eq:noinflcontr}. Contour plots of the density $f(\mathbf{v},t)$ (upper row), and the control $u_1(\mathbf{v},t)$ (bottom row) for time $t=1$ (left), $t=2$ (center) and $t=3$ (right).}
		\label{fig:test3_1}
	\end{figure}
\end{comment}
\begin{comment}
	\begin{figure}[h!]
		\centering\includegraphics[width=0.31\textwidth]{testinfl_T1.eps}
		\quad	\includegraphics[width=0.31\textwidth]{testinfl_T2.eps}\quad \includegraphics[width=0.31\textwidth]{testinfl_T3.eps}
		\\
		\centering
		\includegraphics[width=0.31\textwidth]{infl_contr_T1.eps}\quad\includegraphics[width=0.31\textwidth]{infl_contr_T2.eps}\quad\includegraphics[width=0.31\textwidth]{infl_contr_T3.eps}
		\caption{Case with $s_1(v_1,v_2)$ as in \eqref{eq:inflcontr}. Contour plots of the density $f(\mathbf{v},t)$ (upper row), and the control $u_1(\mathbf{v},t)$ (bottom row) for time $t=1$ (left), $t=2$ (center) and $t=3$ (right).}
		\label{fig:test3_2}
	\end{figure}
\end{comment}

\subsubsection{Bivariate opinion model}\label{sec:dueopi}
\begin{figure}[h!]
	\centering
	\includegraphics[width=0.31\textwidth]{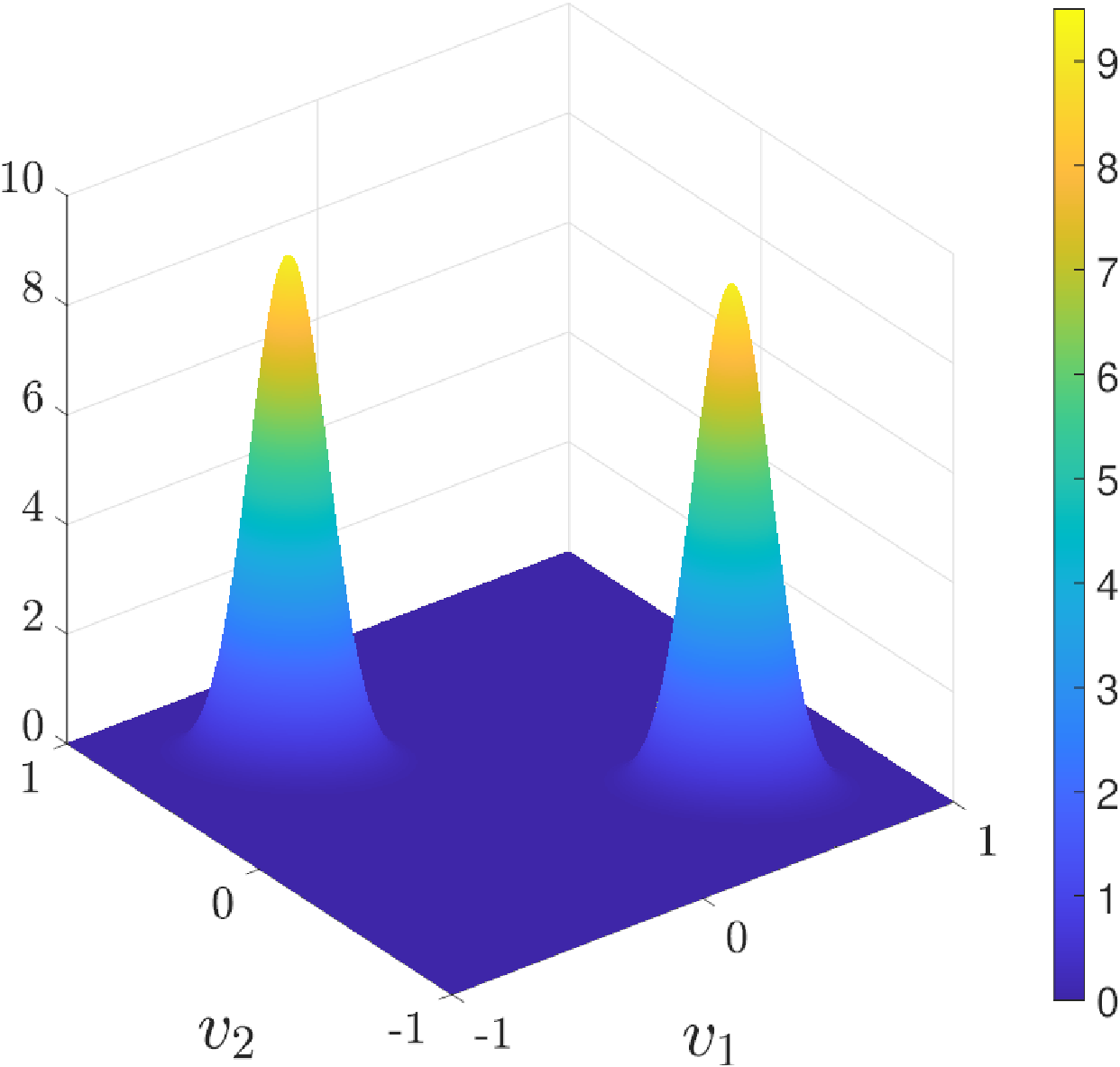}\quad\includegraphics[width=0.31\textwidth]{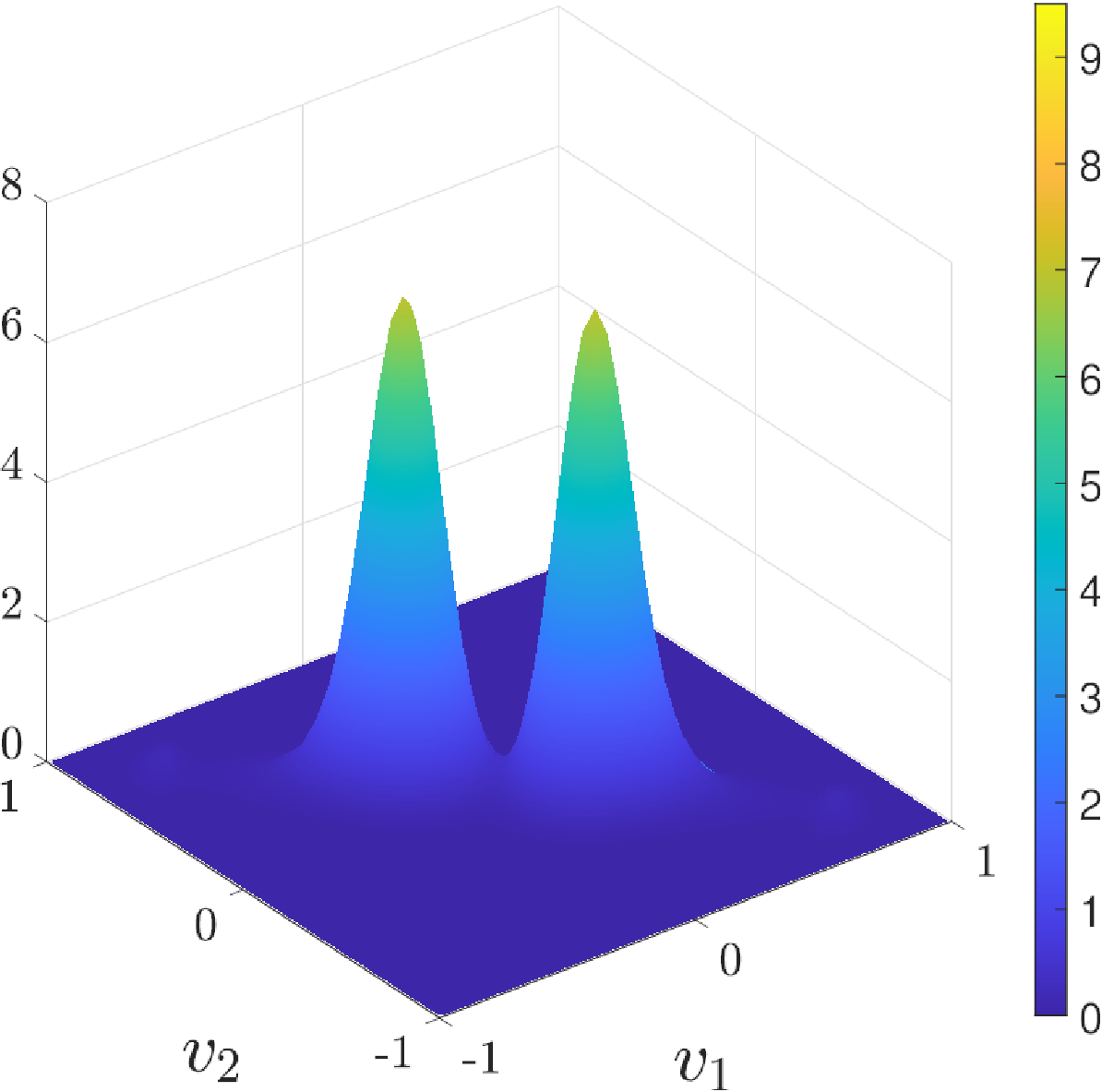}\quad\includegraphics[width=0.31\textwidth]{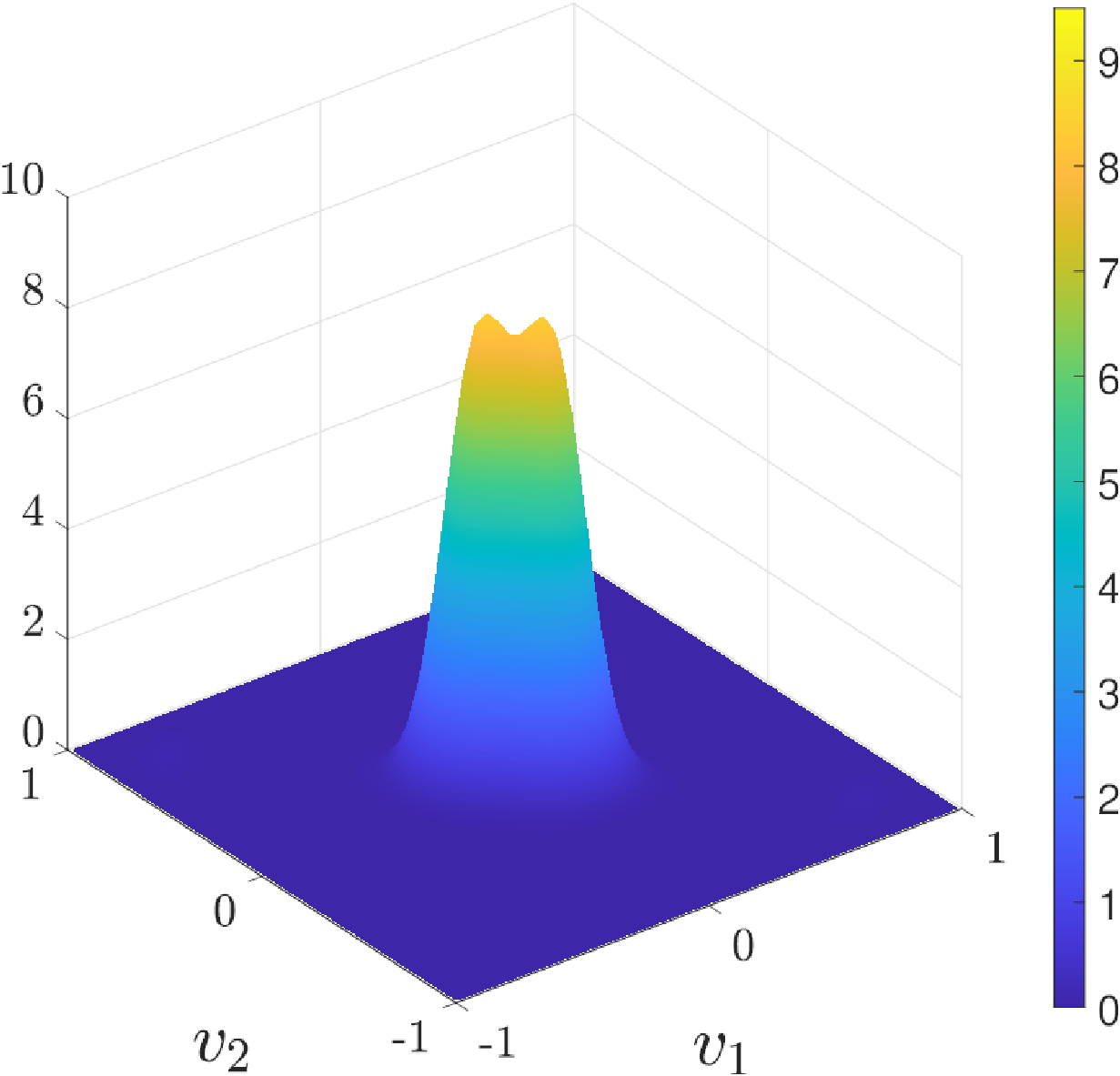}
	\caption{The initial condition $f(\mathbf{v},0)$ (left), the uncontrolled density at time $t=1$ (center) and at time $t=2$ (right).}
	\label{fig:test_2opi2}
\end{figure}

\begin{figure}[h!]
	\centering
	\includegraphics[width=0.3\textwidth]{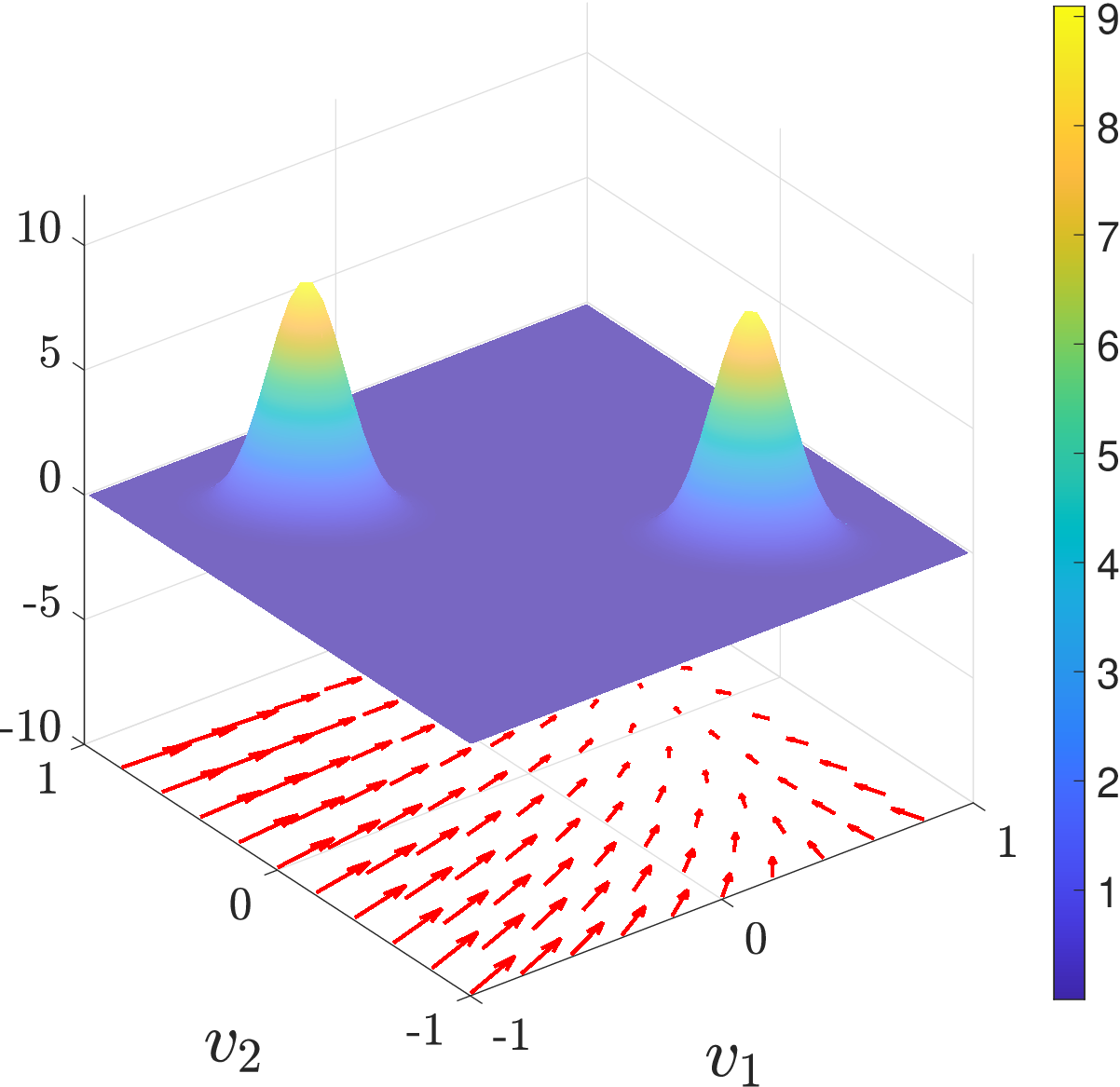}\quad\includegraphics[width=0.3\textwidth]{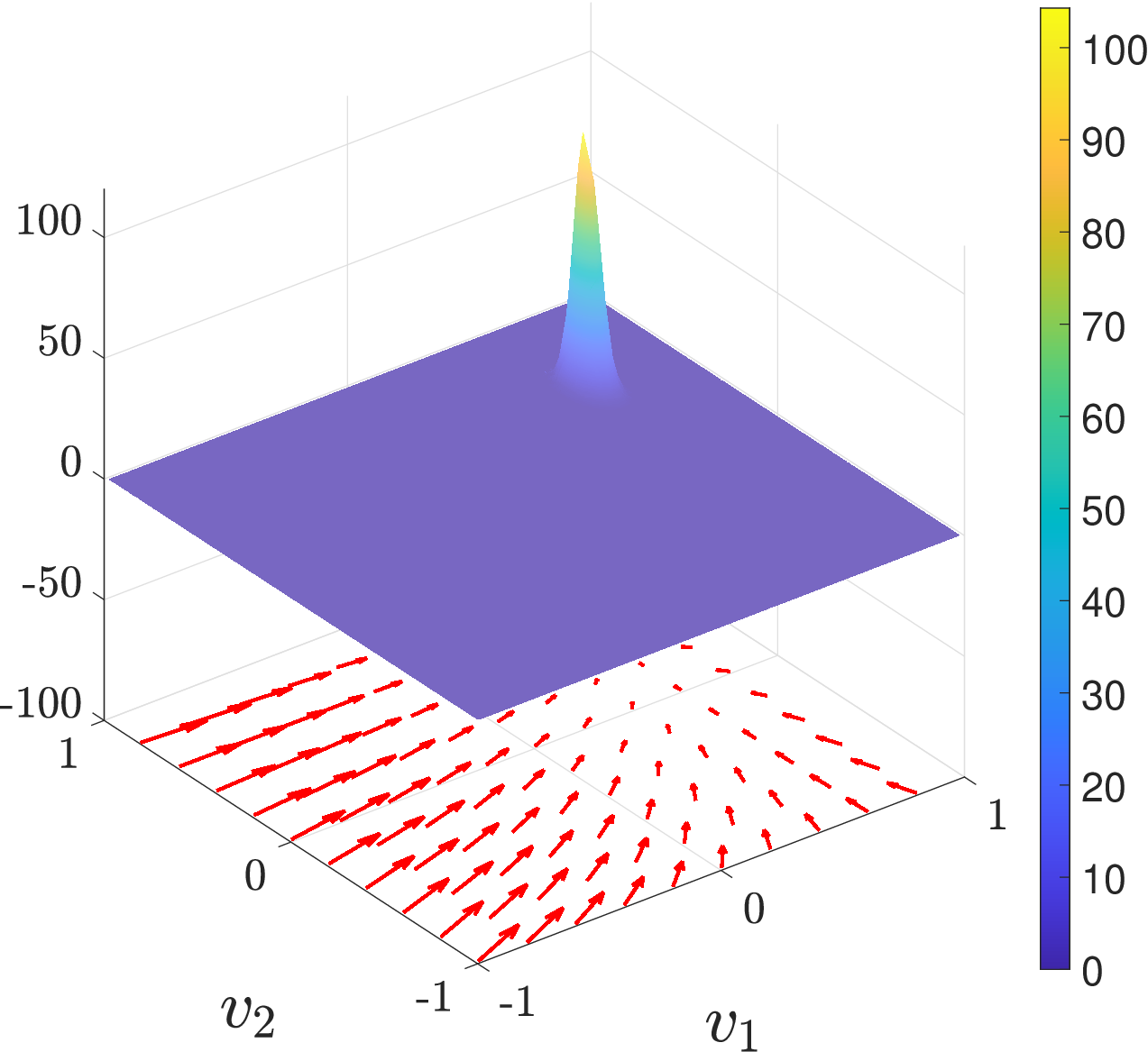}\quad\includegraphics[width=0.3\textwidth]{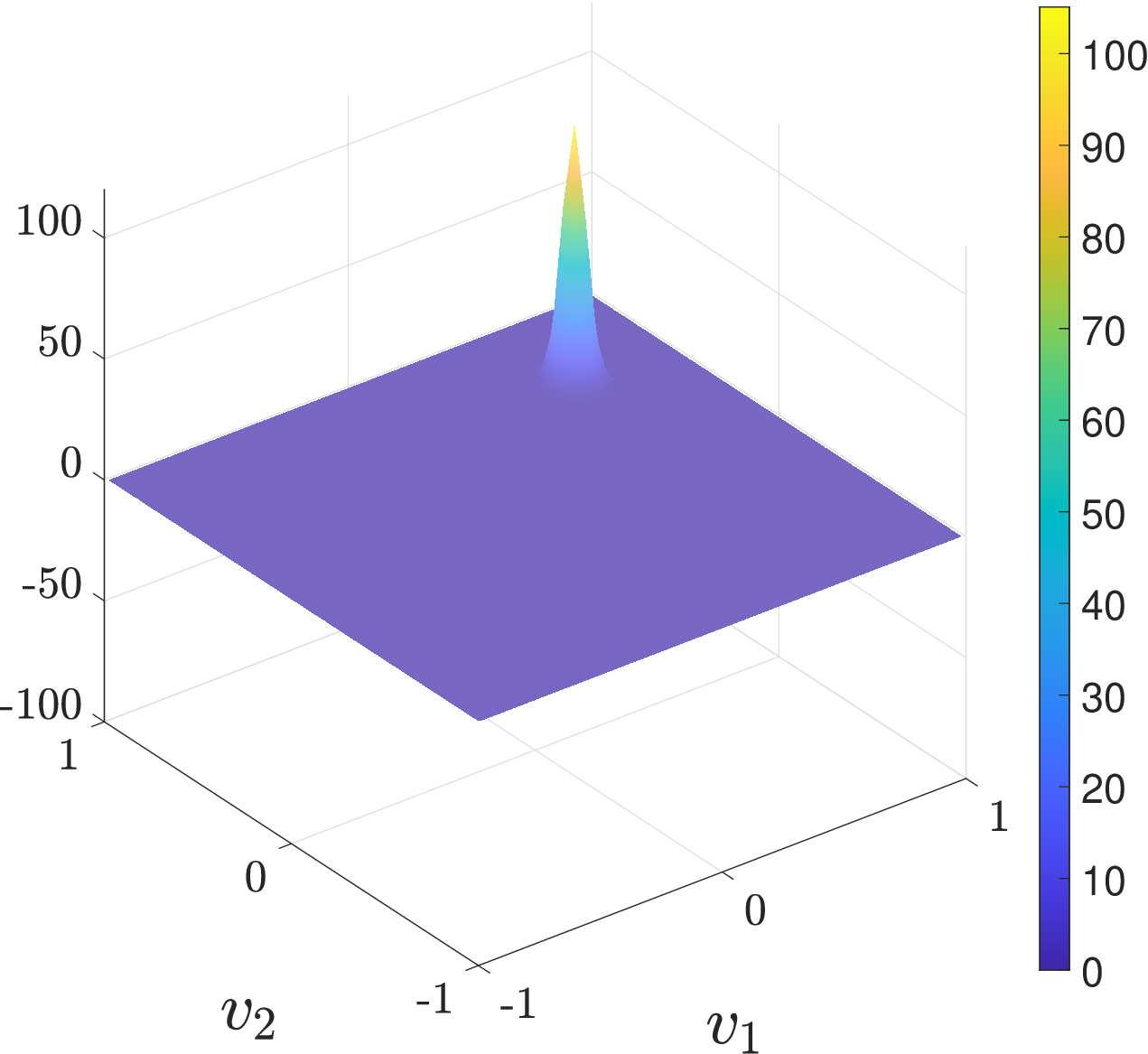}
	\caption{Comparison between the density $f(\mathbf{v},t)$ (as surface plot) and the optimal control $\mathbf{u}^*(\mathbf{v},t)$ (red arrows) at time $t=0$ (left), $t=1$ (center) and $t=2$ (right).}
	\label{fig:test_2opi3}
\end{figure}

\begin{figure}[h!]
	\centering
	\includegraphics[width=0.7\textwidth]{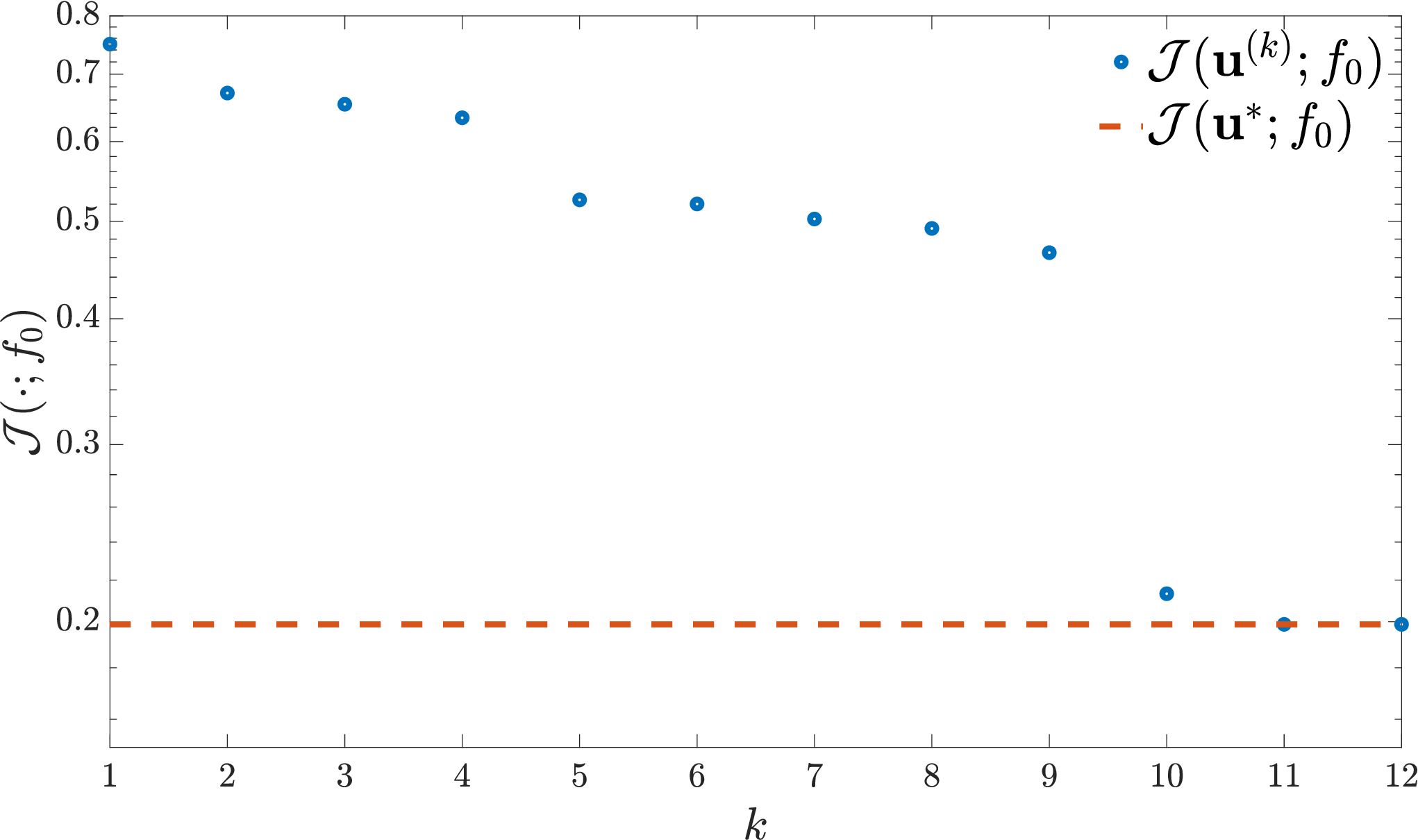}
	\caption{\second{Value of the cost functional $\mathcal{J}(\mathbf{u}^{(k)};f_0)$ at each iteration $k$ of the gradient descent algorithm compared to the final optimal value $\mathcal{J}(\mathbf{u}^*;f_0) = 1.99\times 10^{-1}$.}}
	\label{fig:test_2opi4}
\end{figure}

In this last Section\second{,} we assume that $\Omega_1=\Omega_2=[-1,1]$ and the two variables can be interpreted as opinions. In \eqref{eq:drift_general}\second{,} we set $h_1(v_1,v_2)=h_2(v_1,v_2)=0$, \second{and} 
\begin{equation}
	P_1(v_1,v_2,v_1^*,v_2^*)=P_2(v_1,v_2,v_1^*,v_2^*) = \sqrt{(v_1-v_1^*)^2 + (v_2-v_2^*)^2} \leq \Delta,
\end{equation}
with $\Delta = 1$. We set as diffusion matrix
\[
\mathcal{D}(v_1,v_2)=\begin{bmatrix} \frac{\sigma^2}{2}(1-v_1^2)(1-v_2^2) & 0 \\ 0 & \frac{\nu^2}{2} (1-v_1^2)(1-v_2^2)    
\end{bmatrix},
\]
with $\sigma^2=\nu^2=2\times 10^{-2}$. We fix final time $T = 2$, and in \eqref{eq:cost_j} we set $\overline v_1 = \overline v_2 = 0.6$, and $\gamma = 0.05$. We start from a polarized initial condition given by
\begin{gather*}
	f_0(v_1,v_2) = e^{-\frac{1}{2(1-0.5^2)}\left(\frac{(v_1+0.5)^2}{0.01} + \left(\frac{v_1+0.5}{0.1}\right)\left(\frac{v_2-0.5}{0.1}\right) + \frac{(v_2-0.5)^2}{0.01}\right)} \\+ e^{-\frac{1}{2(1-0.5^2}\left(\frac{(v_1-0.5)^2}{0.01} + \left(\frac{v_1-0.5}{0.1}\right)\left(\frac{v_2+0.5}{0.1}\right) + \frac{(v_2+0.5)^2}{0.01}\right)},
\end{gather*}
meaning that the initial density consists of two separated clusters of opinions.
Figure \ref{fig:test_2opi2} shows a comparison between the initial condition (left), the uncontrolled density at time $t=1$ (center) and the uncontrolled density at the final time $t=2$ (right). We can see how, with our choice of parameters, the uncontrolled setting leads to a consensus around zero both for the first and the second variable. Figure \ref{fig:test_2opi3} shows a comparison between the controlled density and the control at time $t=0$ (left), $t=1$ (center) and $t=2$ (right). We can see how the controlled density concentrates around the target $(\overline v_1,\overline v_2)$ and that the control points towards this point \second{for $t=0$ and $t=1$. The absence of arrows in the rightmost picture is due to the control being equal to zero at the final time $t=2$}. Figure \ref{fig:test_2opi4} shows the value of the cost functional $\mathcal{J}(\mathbf{u}^{(k)};f_0)$ for the various iterations $k$ of the gradient descent method compared to the final computed optimal value $\mathcal{J}(\mathbf{u}^*;f_0)$.

%\begin{figure}[h!]
%	\centering
%	\includegraphics[width=0.3\textwidth]{figs/test4_contr_T0_surf_orange.eps}\quad\includegraphics[width=0.3\textwidth]{figs/test4_contr_T1_surf_orange.eps}\quad\includegraphics[width=0.3\textwidth]{figs/test4_contr_T2_surf_orange.eps}
%	\caption{Comparison between the density $f(\mathbf{v},t)$ (as surface plot) and the optimal control $\mathbf{u}^*(\mathbf{v},t)$ (as quiver plot) at time $t=0$ (left), $t=1$ (center) and $t=2$ (right).}
%	\label{fig:test_2opi3b}
%\end{figure}

\begin{table}[ht]
	\centering\begin{minipage}{.8\textwidth}
		\caption{The table shows the computed optimal value of the cost functional $\mathcal{J}(\mathbf{u}^*;f_0)$, the number of iterations $k$, the CPU time per iteration, the minimum value of the function $f$ and the maximum error in the integral.}
		\label{tab2opi}
		\begin{tabular*}{\textwidth}{@{\extracolsep{\fill}}*{6}{l}@{\extracolsep{\fill}}}
			\toprule%
			$N_v$&$\mathcal{J}(\mathbf{u}^*;f_0)$&$k$&CPU time $(s)$&$f_{min}$ & $E_{int}$\\
			\midrule
			$3$&$2.43\times 10^{-1}$&$10$&$5.32\times 10^{-1}$&$1.66\times 10^{-96}$ & $6.66\times 10^{-16}$\\
			\midrule
			$4$&$1.99\times 10^{-1}$&$11$&$1.90$&$1.64\times 10^{-182}$& $4.44\times 10^{-16}$\\
			\midrule
			$5$&$2.47\times 10^{-1}$&$12$&$1.32\times 10^1$&$8.00\times 10^{-322}$& $1.11\times 10^{-15}$\\
			\bottomrule
	\end{tabular*}\end{minipage}
\end{table}

\section{Conclusions}\label{sec:conclusions}

\second{In this study, we introduced a second-order numerical scheme designed to address optimal control problems governed by nonlinear Fokker–Planck equations, with particular emphasis on applications in socio-economic dynamics and opinion formation. The optimal control problem is formulated by means of Lagrange formalism, resulting in a coupled forward–backward system: a Fokker–Planck equation that models socio-economic phenomena and a Hamilton–Jacobi–Bellman equation governing the adjoint variable. To ensure both high accuracy and stability, we utilize a structure-preserving scheme for the forward equation in conjunction with a semi-Lagrangian approach for the backward equation, thus enabling the adoption of large time steps without the limitations imposed by the CFL condition. The coupling between these two equations is resolved through a reduced-gradient optimization strategy. Numerical experiments confirmed the second-order accuracy of the proposed method and illustrate its efficacy in controlling opinion dynamics across a range of scenarios, including evolving social networks and problems in higher dimensions, e.g. multi-dimensional opinion state.}

\second{The scheme itself is easily implementable also for problems in higher spatial dimensions, however the computational cost increases significantly for each iteration of the gradient method, especially for the semi-Lagrangian step, due to the increasing number of characteristic curves and interpolation needed. It is important to note that, while we can rigorously establish second-order convergence for the discretizations of both the forward and backward equations, we currently lack an analytical proof for the overall accuracy of the method. Consequently the convergence properties of the gradient-based optimization remain supported primarily by empirical evidence.}

\second{Future research directions include the study of Fokker-Planck models with uncertainties, which require the design of robust controls with respect to certain risk measures, and with localized action of the control. A further direction is to investigate the development of numerical schemes able to cope efficiently with the optimal control of partial differential models in presence of multiple scales, \second{these} models are of paramount importance in applications where heterogeneous dynamics are involved, such as opinion dynamics coupled with epidemics, or crowd dynamics with emotional contagion.}

\section*{Acknowledgment}
This work has been written within the activities of the GNCS group of INdAM. GA and EC thanks the support of MUR-PRIN Project 2022 PNRR No. P2022JC95T, “Data-driven discovery and control of multi-scale interacting artificial agent systems”, and GA also of MUR-PRIN Project 2022 No. 2022N9BM3N “Efficient numerical schemes and optimal control methods for time-dependent partial differential equations” both financed by the European Union - Next Generation EU.

%\section*{About References}
%
%Please write at the begining of your {\it bib} file:\\
%\centerline{@settings\{label,options="nameyear"\}}
%so that the references will be sorted in the {\it alphabetical} order according to the surname.
%An example is as follows
%~\cite{bib2,bib1,bib3,bib4,bib5,bib6,bib7,bib8,bib9,bib10,bib11,bib12}.

\bibliographystyle{abbrv}
\bibliography{biblio_opicontrol2}

\end{document}